\documentclass[12pt]{amsart}
\usepackage[all, knot]{xy}
\usepackage{xypic}
\usepackage{epsfig}
\usepackage{xspace}
\usepackage{layout}
\usepackage{latexsym}
\usepackage{amsmath}
\usepackage{amsthm}
\usepackage{amssymb}
\usepackage{amsfonts}
\usepackage{amsxtra}     

\usepackage{verbatim}
\usepackage{longtable}

\input xypic
\input xy
\xyoption{all} \xyoption{arc}

\newcommand{\ie}{\textit{i.e.}\hspace{1.8mm}}

\newcommand{\vcomp}[2]{\genfrac{[}{]}{0pt}{}{#1}{#2}}
\newcommand{\tb}[1]{\phantom{\sum^\Sigma_\Sigma} #1 \phantom{\sum^\Sigma_\Sigma}}
\newcommand{\lr}[1]{\hspace{.5mm}#1\hspace{.5mm}}

\begin{document}

\newtheorem{thm}{Theorem}[section]
\newtheorem{conj}{Conjecture}[section]
\newtheorem{lem}[thm]{Lemma}
\newtheorem{cor}[thm]{Corollary}
\newtheorem{prop}[thm]{Proposition}
\newtheorem{rem}[thm]{Remark}

\theoremstyle{definition}
\newtheorem{defn}[thm]{Definition}
\newtheorem{examp}[thm]{Example}
\newtheorem{notation}[thm]{Notation}
\newtheorem{rmk}[thm]{Remark}

\theoremstyle{remark}

\makeatletter
\renewcommand{\maketag@@@}[1]{\hbox{\m@th\normalsize\normalfont#1}}%
\makeatother
\renewcommand{\theenumi}{\roman{enumi}}

\def\qed {{
   \parfillskip=0pt        
   \widowpenalty=10000     
   \displaywidowpenalty=10000  
   \finalhyphendemerits=0  
  %
   \leavevmode             
   \unskip                 
   \nobreak                
   \hfil                   
   \penalty50              
   \hskip.2em              
   \null                   
   \hfill                  
   $\square$
  %
   \par}}                  

\newenvironment{pf}{{\it Proof:}\quad}{\qed \vskip 12pt}
\date{\today}

\title{Pseudo Algebras and Pseudo Double Categories}
\author[Thomas M. Fiore]{Thomas M. Fiore \\
Department of Mathematics \\
University of Chicago \\
5734 S. University Avenue \\
Chicago, Illinois 60637 \\
USA \\
Phone: 773-702-0088 \\
Fax: 773-702-9787 \\
Email: fiore@math.uchicago.edu}
\address{Department of Mathematics \\ University of Chicago \\ Chicago, IL 60637}
\email{fiore@math.uchicago.edu}
\thanks{The author extends his gratitude to Peter May and Igor Kriz for fruitful
conversations, and also to John Baez, Ronnie Brown, and Michael
Shulman, for some helpful comments on an earlier draft. The author
is indebted to the referee whose careful suggestions greatly
improved this paper.}
\thanks{The author was supported by the Mathematical Sciences Postdoctoral Research
Fellowship of the National Science Foundation. This article was
completed during a visit at the Fields Institute during the Thematic
Program on Geometric Applications of Homotopy Theory, and the author
thanks the Fields Institute for its support.}


\begin{abstract}
As an example of the categorical apparatus of pseudo algebras over
2-theories, we show that pseudo algebras over the 2-theory of
categories can be viewed as pseudo double categories with folding or
as appropriate 2-functors into bicategories. Foldings are equivalent
to connection pairs, and also to thin structures if the vertical and
horizontal morphisms coincide. In a sense, the squares of a double
category with folding are determined in a functorial way by the
2-cells of the horizontal 2-category. As a special case, strict
2-algebras with one object and everything invertible are crossed
modules under a group.
\end{abstract}

\maketitle

\tableofcontents



\section{Introduction}
Recent years have seen widespread applications of {\it
categorification}. The term categorification refers to a process of
turning algebraic notions on sets into algebraic notions on
categories as explained in \cite{baez}. Generally speaking, one
takes a set-based algebraic notion, then replaces sets by
categories, functions by functors, and equations by natural
isomorphisms which satisfy certain coherence diagrams.

For example, a monoid (group without inverses) is a set-based
algebraic concept. Its categorified notion is a monoidal category,
that is a category $M$ equipped with a functor $\xymatrix@1{\otimes:
M \times M \ar[r] & M}$ and a unit which satisfy the monoid axioms
up to coherence isos. These coherence isos must satisfy certain
coherence conditions, such as the familiar pentagon diagram. The
commutativity of these diagrams in turn implies that all diagrams in
a certain class commute, as Mac Lane proved in \cite{maclane1}.  A
familiar example of a monoidal category is the category of complex
vector spaces under the operation of tensor product with unit
$\mathbb{C}$.

Another example of categorification is the notion of a {\it
bicategory}, which is a categorification of the algebraic concept of
category. In a bicategory the hom-sets are categories and
composition is a functor. Composition is unital and associative up
to coherence isomorphisms which satisfy coherence diagrams like
those of a monoidal category. This similarity is not a coincidence:
one-object bicategories are monoidal categories in the same way that
one-object categories are monoids. A familiar example of a
bicategory consists of rings, bimodules, and bimodule morphisms.
Bicategories were introduced in the 1960's in \cite{benabou1},
\cite{benabou}, \cite{ehresmann}, and \cite{ehresmann2}. Since then,
they (and their variants) have appeared in diverse areas, such as
homotopy theory and high energy physics.

However, the question arises: what exactly does one mean by
``coherence isos satisfying certain coherence diagrams''? Which
coherence isos and which coherence diagrams does one require? This
question already suggests that there may be more than one way to
categorify a given concept, such as category. Indeed, there already
are a dozen or so different definitions of weak $n$-category, many
of which are described in \cite{cheng1} and \cite{leinster1}.

Lawvere theories and 2-theories provide one answer to this question.
Lawvere theories, first introduced in \cite{lawvere}, abstractly
encode algebraic structure. For most familiar algebraic structures
there is a Lawvere theory. For example, there is a Lawvere theory of
monoids, and algebras over this theory are precisely the monoids. A
Lawvere theory $T$ is simply a category whose objects are $0,1,2,
\dots$ such that $n$ is the product of $n$ copies of $1$ with
specified projection maps. If $T$ is the theory which encapsulates a
certain algebraic structure, then a set $X$ with that algebraic
structure is an {\it algebra over the theory $T$}. This means that
$X$ is equipped with a morphism $\xymatrix@1{\Phi:T \ar[r] &
End(X)}$ of theories from $T$ to the endomorphism theory on $X$. To
each abstract word $\xymatrix@1{w:n \ar[r] & 1}$ a morphism assigns
a function $\xymatrix@1{\Phi(w):X^n \ar[r] & X}$ in a uniform way.

Similarly a category $X$ is a {\it pseudo algebra over a theory $T$}
if it is equipped with a pseudo morphism of theories
$\xymatrix@1{\Phi:T \ar[r] & End(X)}$. To each abstract word
$\xymatrix@1{w:n \ar[r] & 1}$ a pseudo morphism assigns a functor
$\xymatrix@1{\Phi(w):X^n \ar[r] & X}$. Additionally, for each
operation of theories, there is a coherence isomorphism and for each
relation of theories, there is a coherence diagram which these
coherence isomorphisms must satisfy. This is a well-defined
procedure which specifies exactly which coherence isomorphisms and
coherence diagrams are appropriate, no matter if one is interested
in monoids, semi-rings, rings, etc. A {\it pseudo monoid, pseudo
semi-ring,} or {\it pseudo ring} is simply a pseudo algebra over the
appropriate theory. There is a systematic way to leave out some
coherence diagrams to encompass more examples \cite{fiorekrizhu}.

However, Lawvere theories only axiomatize algebraic structures on a
single set. There is no Lawvere theory of categories, since a
category consists of two sets with composition defined in terms of
pullback. For algebraic structures on several sets, one can use
limit theories, sketches, and multi-sorted theories as in
\cite{adamekrosicky1994}, \cite{bergnermultisorted}, or
\cite{boardmanvogt1973}, or schemes of operators as in \cite{higgins1963}.
In this paper we consider categories as
algebras over a {\it 2-theory}. This adds a new ingredient to
categorification that we do not see in the one-object case of
Lawvere theories. For example, pseudo algebras over the 2-theory of
categories have an object category $I$ instead of an object set, as
we shall see.

This version of categorification in terms of pseudo algebras over
2-theories was introduced in \cite{hu}, and further developed in
\cite{fiore1} and \cite{hu1}, to give a completely rigorous approach
to conformal field theory with $n$-dimensional modular functor.
Pseudo algebras over the {\it 2-theory of commutative monoids with
cancellation} make the symmetric approach to conformal field theory
outlined in \cite{segal1} rigorous. The notion of 2-theory was the
main ingredient for a well-defined procedure of passing from a
strict algebraic structure on a {\it family} of sets to a pseudo
algebraic structure on a {\it family} of categories, such as the
pseudo algebraic structure of disjoint union and gluing on the class
of worldsheets (rigged surfaces). This procedure gave a well-defined
machine for generating the coherence isos and coherence diagrams
that were missing from conformal field theory until that point.
Already in 1991, Mac Lane suggested a study of coherence in the
context of conformal field theory in \cite{maclane1991}. The
foundations of pseudo algebras over theories and 2-theories were
written in \cite{fiore1}, as well as theorems relevant for
application to conformal field theory. Among other things, it was
shown that 2-categories of pseudo algebras admit pseudo limits and
bicolimits, and forgetful 2-functors of pseudo algebras admit left
biadjoints.

In the present article we apply this version of categorification to
the fundamental algebraic structure of category and compare the
resulting concept of pseudo category to weak double categories and
also pseudo functors $\xymatrix@1{I \ar[r] & \mathcal{C}}$. One
might expect that a pseudo category would neatly fit into one of the
two prevailing approaches to categorification: enrichment and
internalization. This however is not true, a pseudo category is
neither a bicategory, nor a weak double category. Instead we arrive
at an intermediate notion: a pseudo category can be 2-equivalently
described as a weak double category with weak folding or as a
bicategory equipped with a pseudo functor from a 1-category. Our
pseudo categories are slightly different from the pseudo categories
in \cite{martinspseudocategories}, so we will call them pseudo
$I$-categories instead.

We first treat the categorified strict case by reviewing strict
categories and double categories in Section \ref{strictcategories}
and Section \ref{doublesection}, and prove the strict versions of
our desired result in Theorems \ref{YZ}, \ref{XY}, and \ref{XZ}.
Foldings, used in \cite{brownmosa99}, are introduced to facilitate
the 2-equivalence of strict 2-algebras over the 2-theory of
categories with underlying category $I$ ($I$-categories for short) and certain double categories.
It turns out that foldings, which have Ehresmann's quintets as their
motivating example, are equivalent to Brown and Spencer's connection
pairs, and also thin structures in the edge-symmetric case, as
recounted in Theorem \ref{connection=horizontalization} (Lemmas
\ref{connectionpairtohorizontalization}-\ref{lastpiece}) and
Corollary \ref{folding=thin}. In light of this, Theorem
\ref{YZ} is an $I$-category analogue of the equivalence in \cite{brownmosa99} and
\cite{spencer77} between small 2-categories and edge-symmetric
double categories with thin structure.

In the case of one object with everything invertible, strict
2-algebras (not necessarily edge symmetric) are equivalent to
crossed modules under a group as in Theorem \ref{XYZinv} and Theorem
\ref{WZ}. We generalize Brown and Spencer's equivalence in
\cite{brownspencer76} between edge-symmetric double groups with
connection pair and crossed modules. In Theorem
\ref{specialextension} we prove that double groups (not necessarily
edge symmetric) with folding are 2-equivalent to crossed modules
under groups. The paper \cite{brownmackenzie} contains a substantial
generalization of \cite{brownspencer76} by giving an equivalence of
`core diagrams' to double groupoids with certain filling conditions.
Double groupoids have recently found application in the theory of
weak Hopf algebras in \cite{andruskiewitschnatalequantum} and
\cite{andruskiewitschnataletensor}.

The pseudo double categories of \cite{grandisdouble1} are reviewed
in Section \ref{pseudodoublesection}. We finally prove in Theorem
\ref{XYpseudo} and Theorem \ref{XZpseudo}, under the assumption of
strict units, the 2-equivalence of pseudo algebras over the 2-theory
of categories (pseudo $I$-categories for short), pseudo double
categories with folding, and strict 2-functors\footnote{The term
``2-functor'' means strict 2-functor in this paper. Sometimes we
include the word ``strict'' for emphasis. When we mean pseudo
functor (homomorphism of bicategories), we say so.} from a groupoid
into a bicategory. The latter two 2-categories remain 2-equivalent
even if $I$ is merely a category.

Theorem \ref{XZ} and Theorem \ref{XZpseudo} may also be considered a
special case of Theorem 6.5 of the comparison article \cite{fiore3}.
That article relates the commutative-monoid-with-cancellation
approach to conformal field theory in \cite{hu} (outlined in
\cite{segal1} in terms of trace) to the cobordism approach.

\section{Strict $I$-Categories}\label{strictcategories}

A category consists of a family of sets with an algebraic structure.
Namely, if $\mathcal{C}$ is a category with object set $I$, then the
associated family of sets $X_{A,B}:=Hom_\mathcal{C}(A,B)$ is
parametrized by $I^2$. On this family of sets, we have the algebraic
structure of composition and identity. Thus we can view a category
as a functor $\xymatrix@1{X:I^2 \ar[r] & Sets}$ with certain
algebraic operations, where $I^2$ is considered as a discrete
category.

From this point of view, a category is an algebra $\xymatrix@1{X:I^2
\ar[r] & Sets}$ over the {\it 2-theory of categories}. This is the
notion that we categorify. In this article we do not write down the
2-theory of categories, since it suffices to directly define
2-algebras and pseudo algebras over this 2-theory. The operations
are given in terms of the generating words $\circ$ and $\eta$ rather
than abstract operations of 2-theories. The underlying theory of the
2-theory of categories is the theory of sets. We take the following
description as a definition, and do not need the notion of 2-theory
anywhere in this paper. For a development of 2-theories and their
algebras, see the original paper \cite{hu}, or the papers
\cite{fiore1} and \cite{hu1}.

\begin{defn} \label{algebradefinition}
A {\it strict 2-algebra}\footnote{The distinction in this paper
between strict 2-algebra and pseudo algebra agrees with usual
2-terminology and pseudo terminology. For example a strict 2-functor
is a morphism of strict 2-categories. A pseudo functor is on the
same level as a 2-functor, except a pseudo functor preserves
composition and unit only up to coherent 2-cell isomorphisms. The
notion of pseudo 2-algebra over a theory in \cite{fiore3} is
distinct from a pseudo algebra over a theory. It should also be
noted that a 2-theory is not a theory enriched in categories, nor
any sort of weakened theory.} {\it over the 2-theory of categories with underlying category $I$},
called {\it $I$-category} for short, consists of a category $I$ and a
strict 2-functor $\xymatrix@1{X:I^2 \ar[r] & Cat}$ with strictly
2-natural functors
$$\xymatrix@1{X_{B,C} \times X_{A,B} \ar[r]^-{\circ} & X_{A,C}}$$
$$\xymatrix@1{\ast \ar[r]^{\eta_B} & X_{B,B} }$$
for all $A,B,C \in I$. These functors satisfy the following
relations.

\begin{enumerate}
\item
The composition $\circ$ is {\it associative}.
\begin{small}
$$\xymatrix@C=5pc{(X_{C,D} \times X_{B,C}) \times X_{A,B} \ar[r]^-{\circ \times 1_{X_{A,B}}}
\ar[dd]_{\cong} & X_{B,D} \times
X_{A,B} \ar[dr]^{\circ} & \\
& &  X_{A,D} \\ X_{C,D} \times (X_{B,C} \times X_{A,B})
\ar[r]_-{1_{X_{C,D}} \times \circ} & X_{C,D} \times X_{A,C}
\ar[ur]_{\circ} & }$$
\end{small}
\item
For each $B \in I$, the operation $\eta_B$ is an {\it identity} for
$\circ$.

\begingroup
\vspace{-2\abovedisplayskip} \Small
$$\begin{array}{cl}
\xymatrix@C=4pc@R=3pc{\ast \times X_{A,B} \ar[r]^-{\eta_B \times
1_{X_{A,B}}} \ar[dr]_{pr_2} & X_{B,B} \times X_{A,B} \ar[d]^-{\circ}
\\ & X_{A,B} } & \xymatrix@C=4pc@R=3pc{X_{B,C} \times \ast
\ar[r]^-{1_{X_{B,C}} \times \eta_B} \ar[dr]_{pr_1} & X_{B,C} \times
X_{B,B} \ar[d]^-{\circ} \\ & X_{B,C} }
\end{array}$$
\endgroup
\noindent
\end{enumerate}
We denote the value of $\eta_B$ on the unique object and morphism of
the terminal category by $1_B$ and $i_{1_B}$ respectively. We denote
the identity morphism on an object $f$ in the category $X_{A,B}$ by
$i_f$.
\end{defn}


The term {\it $I$-category} is an abbreviation of
{\it strict 2-algebra over the 2-theory of categories with underlying category $I$}. The strict
morphisms of $I$-categories and their 2-cells below are the strict
morphisms and 2-cells in the 2-category of strict 2-algebras over
the 2-theory of categories with the same underlying groupoid $I$ as
in \cite{fiore1}, \cite{hu}, and \cite{hu1}.

The term {\it $I$-category} agrees with existing usage of the
term {\it $O$-category} to mean a category with the object set $O$.
Indeed, if $I$ is a discrete category (\ie a set) and $X$ takes
values in $Sets$, then an $I$-category is precisely an ordinary
category with object set $I$.  More generally for groupoids $I$, we
will see that $I$-categories are ``categories with object groupoid
$I$'' in a precise sense.

\begin{defn}
A {\it strict morphism $\xymatrix@1{F: X \ar[r] & Y}$ of
$I$-categories} is a strict 2-natural transformation
$\xymatrix@1{F:X \ar@{=>}[r] & Y}$ which preserves composition and
identity strictly.
\end{defn}

\begin{defn}
A {\it 2-cell $\xymatrix@1{\sigma:F \ar@{=>}[r] & G}$} in the
2-category of $I$-categories is a modification $\xymatrix@1{\sigma:F
\ar@{~>}[r] & G}$ compatible with composition and identity. More
specifically, a 2-cell $\sigma$ consists of natural transformations
$\xymatrix@1{\sigma_{A,B}:F_{A,B} \ar@{=>}[r] & G_{A,B}}$ for all
$A,B \in I$ such that
$$Y_{j,k}(\sigma_{A,B}^f)=\sigma_{C,D}^{X_{j,k}(f)}$$
$$\sigma^g_{B,C} \circ \sigma^f_{A,B}=\sigma_{A,C}^{g \circ f}$$
$$\sigma_{A,A}^{1_A}=i_{1_A}$$ for all $\xymatrix@1{(j,k):(A,B) \ar[r] &
(C,D)}$ in $I^2,$ $f \in X_{A,B},$ and $g \in X_{B,C}$. Here $\circ$
denotes the composition functor of the strict 2-algebra, {\it not}
the composition in the categories $X_{A,B}$.
\end{defn}

\begin{lem}
If $\xymatrix@1{X:I^2 \ar[r] & Cat}$ is an $I$-category and $I$ is a
discrete category, then $X$ is a strict 2-category with object set
$I$. A morphism between two such $I$-categories is simply a strict
2-functor which is the identity on objects. A 2-cell
$\xymatrix@1{\sigma:F \ar@{=>}[r] & G}$ is an oplax natural
transformation with identity components. If additionally $X$ takes
values in $Sets$, then $X$ is simply a category with object set $I$
in the usual sense. A morphism between such $I$-categories is simply
a functor which is the identity on objects. There is at most a
trivial 2-cell between any two such morphisms.
\end{lem}
\begin{pf}
This is just a matter of definitions. The category of morphisms from
$A$ to $B$ is $X_{A,B}$.
\end{pf}

An $I$-category is {\it not} an internal category in $Cat$, {\it
nor} a $Cat$-enriched category, since we have taken as our starting
point a different description of category. More specifically, if one
takes as a starting point the definition of category as an object
set $C_0$ and an arrow set $C_1$ along with four maps defining
source, target, identity, and composition satisfying the relevant
axioms, then one indeed arrives at the notion of internal category
in $Cat$ as described on pages 267-270 of \cite{maclane3}. An
internal category in $Cat$ is the same as a double category, which
is described in an elementary way in the next section. The choice of
starting point is crucial for higher-dimensional category theory. As
seen in \cite{cheng1}, equivalent definitions of category lead to
quite different notions of higher category. We will see that the
2-cells of $I$-categories correspond to certain vertical natural
transformations between double functors.

The notion of $I$-category lies between the notions of internal
category in $Cat$ and $Cat$-enriched category, so how far away is an
$I$-category from a 2-category? The following Lemma shows how to
associate to an $I$-category a strict 2-functor $\xymatrix@1{P:I
\ar[r] & \mathcal{C}}$. More importantly, in the presence of the
other 2-cell axioms, we obtain a simplification of the requirement
that 2-cells $\sigma$ be modifications in terms of compatibility
with $P$. The 2-equivalence of $I$-categories to such strict
2-functors is Theorem \ref{XZ}. We will also apply the following
Lemma in the comparison with double categories with folding in
Theorem \ref{XY}.

\begin{lem} \label{holonomyconstruction}
Suppose $I$ is a groupoid, $X$ and $X'$ are strict $I$-categories,
and $\xymatrix@1{F,G:X \ar[r] & X'}$ are strict morphisms. We
associate to $X$ a 2-category $\mathcal{C}$ with $Obj \hspace{1mm}
\mathcal{C}:=Obj \hspace{1mm} I$ and
$Mor_{\mathcal{C}}(A,B):=X_{A,B}$. We denote the identity on $A$ in
the category $I$ by $1^v_A$ while we denote the identity on $A$ in
$\mathcal{C}$ by $1^h_A$. The identity 2-cell in $\mathcal{C}$ on a
morphism $f$ is $i_f$. Let $\xymatrix@1{P:I \ar[r] & \mathcal{C}}$
be the strict 2-functor which is the identity on objects and
$$P(j):=X_{j^{-1},1^v_C}(1^h_C)=X_{1^v_A,j}(1^h_A)$$ for morphisms
$j \in I(A,C)$. Let $\mathcal{C}'$ and $\xymatrix@1{P':I \ar[r] &
\mathcal{C}'}$ be the 2-category and strict 2-functor associated
analogously to $X'$. Suppose further we have for each $A,B \in I$ a
natural transformation $\xymatrix@1{\sigma_{A,B}:F_{A,B} \ar@{=>}[r]
& G_{A,B}}$ such that
$$\sigma^g_{B,C} \circ \sigma^f_{A,B}=\sigma_{A,C}^{g \circ f}$$
$$\sigma_{A,A}^{1_A^h}=i_{1_A^h}$$
for all $f \in X_{A,B}$ and $g \in X_{B,C}$. Then the following are
equivalent.
\begin{enumerate}
\item
For all $\xymatrix@1{(j,k):(A,B) \ar[r] & (C,D)}$ in $I^2$ and all
$f \in X_{A,B}$ we have
$Y_{j,k}(\sigma_{A,B}^f)=\sigma_{C,D}^{X_{j,k}(f)}$.
\item
For all $\xymatrix@1{j:A \ar[r] & C}$ in $I$ we have
$\sigma_{A,C}^{P(j)}=i_{P'(j)}$.
\end{enumerate}
\end{lem}
\begin{pf}
The naturality of the identity implies
$$X_{j^{-1},1^v_C}(1^h_C)=X_{1^v_A,j}(1^h_A).$$ The map $P$ preserves
compositions $\xymatrix@1{A \ar[r]^{j} & C \ar[r]^{k} & E}$ in $I$
because
$$
\aligned
P(k \circ j) &= X_{1^v_A,k \circ j}(1_A^h) \\
&=X_{1^v_A,k}(X_{1^v_A,j}(1_A^h)) \\
&=X_{1^v_A,k}(P(j)) \\
&=X_{1^v_A,k}(1_C^h \circ P(j)) \\
&=(X_{1^v_C,k}(1_C^h)) \circ P(j) \\
&=P(k) \circ P(j)
\endaligned
$$
by the naturality diagram below.
$$\xymatrix@R=3pc@C=3pc{X_{C,C} \times X_{A,C} \ar[r]^-{\circ}
\ar[d]_{X_{1^v_C,k} \, \times \, X_{1^v_A,1^v_C}} & X_{A,C}
\ar[d]^{X_{1^v_A,k}}
\\ X_{C,E} \times X_{A,C} \ar[r]_-{\circ} & X_{A,E}}$$
It is clear that $P(1^v_A)=X_{1^v_A,1^v_A}(1^h_A)=1^h_A$, so we
indeed have a 2-functor $P$.

For $f \in X_{A,B}$ and $j\in I(A,C)$, note that $f \circ
P(j^{-1})=X_{j,1^v_B}(f)$ by the naturality diagram
$$
\xymatrix@R=3pc@C=3pc{ X_{A,B} \times X_{C,A} \ar[r]^-{\circ}
\ar[d]_{X_{j,1_B^v} \times \ X_{1_C^v,j}} & X_{C,B}
\ar[d]^{X_{1_C^v,1_B^v}}
\\ X_{C,B} \times X_{C,C} \ar[r]_-{\circ} & X_{C,B}}
$$
and similarly $P(k) \circ f = X_{1^v_A,k}(f)$ for $k\in I(B,D)$.
Similar statements hold for 2-cells of $\mathcal{C}$. Thus
$$X_{j,k}(f)=P(k)\circ f \circ P(j^{-1})$$
$$X_{j,k}(\alpha)=i_{P(k)} \circ \alpha \circ i_{P(j^{-1})}.$$
We use $\circ$ to denote the horizontal composition of 2-cells in a
2-category, in addition to the composition of morphisms.

Let $\sigma_{A,B}$ be a natural transformation for each $A,B \in I$
such that $\sigma$ is compatible with composition and identity.
Suppose $\sigma$ satisfies (i). Then
$$\aligned
\sigma_{A,C}^{P(j)} & = \sigma_{A,C}^{X_{1^v_A,j}(1^h_A)} \\
& = Y_{1_A^v,j}(\sigma^{1^h_A}_{A,A}) \\
& = Y_{1_A^v,j}(i_{1^h_A}) \\
& = i_{P'(j)}
\endaligned$$
and (ii) holds.

Suppose $\sigma$ satisfies (ii). Then
$$\aligned
Y_{j,k}(\sigma^f_{A,B}) & = i_{P'(k)} \circ \sigma^f_{A,B} \circ i_{P'(j^{-1})} \\
&= \sigma^{P(k)}_{B,D} \circ \sigma^f_{A,B} \circ \sigma^{P(j^{-1})}_{C,A} \\
&= \sigma^{P(k)\circ f \circ P(j^{-1})}_{C,D} \\ &=
\sigma^{X_{j,k}(f)}_{C,D}
\endaligned$$
and (i) holds.
\end{pf}

\section{Double Categories with Folding} \label{doublesection}

Ehresmann introduced double categories in \cite{ehresmann} and
\cite{ehresmann2}. After a long gestation period, a full theory of
double categories is beginning to emerge. Classics in the subject
include \cite{ehresmannone}, \cite{brownspencer76},
\cite{brownspencer74},
\cite{ehresmanntwo}-\cite{ehresmann2}\nocite{ehresmannthree}\nocite{ehresmannfour},
and \cite{lodayfinitelymany}.  For recent work on double categories
and related topics, see \cite{alglbrownsteiner2002},
\cite{brown99}-\nocite{browngilbert1989}\nocite{brownicen2003}\cite{brownmosa99},
\cite{dawsonparepronk2003adjoining}-\cite{dawsonpare2002free}\nocite{dawsonparepronkundecidable2003}\nocite{dawsonparepronkfree2004}\nocite{dawsonparepronkpathology2006}\nocite{dawsonparepronkspan2006}\nocite{dawsonpare1993},
\cite{garner2005}-\nocite{garnerthesis}\nocite{garnerpseudodistributive2005}\nocite{grandisdouble1}\cite{grandisdouble2},
\cite{leinsteroperads2004}, and
\cite{martins2004}-\nocite{martinspseudocategories}\cite{mortondouble}.

We recall double categories and foldings, as well as their morphisms
and transformations. Foldings allow us to compare double categories
with $I$-categories in the next section. In Theorem
\ref{connection=horizontalization} we show that foldings are
equivalent to connection pairs, as a corollary they are also
equivalent to thin structures in the edge-symmetric case.

\begin{defn}
A {\it double category} $\mathbb{D}=(\mathbb{D}_0,\mathbb{D}_1)$ is
a category object in the category of small categories. This means
$\mathbb{D}_0$ and $\mathbb{D}_1$ are categories equipped with
functors $$\xymatrix@C=3pc{\mathbb{D}_1 \times_{\mathbb{D}_0}
\mathbb{D}_1  \ar[r] & \mathbb{D}_1 \ar@/^1pc/[r]^s \ar@/_1pc/[r]_t
& \ar[l]|{\lr{\eta}} \mathbb{D}_0 }$$ that satisfy the usual axioms
of a category. We call the objects and morphisms of $\mathbb{D}_0$
the {\it objects} and {\it vertical morphisms} of $\mathbb{D}$, and
we call the objects and morphisms of $\mathbb{D}_1$ the {\it
horizontal morphisms} and {\it squares} of $\mathbb{D}$.
\end{defn}

We can expand this definition as in \cite{kelly}. A {\it double
category $\mathbb{D}$} consists of a set of {\it objects}, a set of
{\it horizontal morphisms}, a set of {\it vertical morphisms}, and a
set of {\it squares} equipped with various sources, targets, and
associative and unital compostions as follows. Objects are denoted
with capital Latin letters $A,B,\dots,$ horizontal morphisms are
denoted with lower-case Latin letters $f,g, \dots,$ vertical
morphisms are denoted with lower-case Latin letters $j,k,\dots,$ and
squares are denoted with lower-case Greek letters
$\alpha,\beta,\dots$ with source and target as indicated below.
\begin{equation} \label{sourcetarget}
\xymatrix{A \ar[r]^f &  B & A \ar[d]_j
& A \ar[r]^f \ar[d]_j \ar@{}[dr]|\alpha & B \ar[d]^k \\
& & C & C \ar[r]_g & D}
\end{equation}
In particular, $\alpha$ has vertical source and target $f$ and $g$,
and horizontal source and target $j$ and $k$ respectively. The
objects and vertical morphisms form a category with composition
denoted $$j_2 \circ j_1 = \begin{bmatrix} j_1 \\ j_2 \end{bmatrix}$$
and identities denoted $1^v_A$. The objects and horizontal morphisms
also form a category, with composition denoted $$f_2 \circ f_1 = [
f_1 \ f_2  ]$$ and identities $1^h_A$. The vertical morphisms and
squares form a category under horizontal composition of squares,
with horizontal identity squares denoted
\begin{equation} \label{horizontalidentity}
\xymatrix{A \ar[r]^{1^h_A} \ar[d]_j \ar@{}[dr]|{i_j^h} & A \ar[d]^j  \\
C \ar[r]_{1^h_C} & C.}
\end{equation}
If $\alpha$ and $\beta$ are horizontally composable squares, then
their composition is denoted $$[ \alpha \ \beta  ].$$ The horizontal
morphisms and squares form a category under vertical composition of
squares, with vertical identity squares denoted
\begin{equation} \label{verticalidentity}
\xymatrix{A \ar[r]^f \ar@{}[dr]|{i_f^v} \ar[d]_{1^v_A} & B
\ar[d]^{1^v_B} \\ A \ar[r]_f & B.}
\end{equation}
If $\alpha$ and $\beta$ are vertically composable squares, then
their composition is denoted $$\begin{bmatrix} \alpha \\ \beta
\end{bmatrix}.$$
The identity squares are compatible with horizontal and vertical
composition.
$$
\begin{array}{c}
\begin{bmatrix}i_{f_1}^v & i_{f_2}^v \end{bmatrix} = i_{[f_1 \ f_2]}^v\end{array}
\hspace{1in}
\begin{array}{c}
\begin{bmatrix} i_{j_1}^h \vspace{2mm} \\ i_{j_2}^h \end{bmatrix}=\ i_{\vcomp{j_1}{j_2 \vspace
{1mm}}}^h.
\end{array}$$
Lastly, the {\it interchange law} holds, \ie in the situation
\begin{equation} \label{interchange1}
\xymatrix{\ar[r] \ar[d] \ar@{}[dr]|{\alpha} & \ar[r] \ar[d] \ar@{}[dr]|{\beta} & \ar[d]  \\
\ar[r] \ar[d] \ar@{}[dr]|{\gamma} & \ar[r] \ar[d]
\ar@{}[dr]|{\delta} & \ar[d]
\\\ar[r] & \ar[r]  &  }
\end{equation}
we have
$$\begin{bmatrix}
\begin{bmatrix}
\alpha & \beta
\end{bmatrix} \vspace{1mm} \\
\begin{bmatrix}
\gamma & \delta
\end{bmatrix}
\end{bmatrix}=\begin{bmatrix}
\begin{bmatrix} \alpha \\ \gamma \end{bmatrix} &
\begin{bmatrix} \beta \\ \delta \end{bmatrix}
\end{bmatrix}$$
and this composition is denoted
\begin{equation} \label{interchange2}
\begin{bmatrix} \alpha & \beta \\ \gamma & \delta
\end{bmatrix}.
\end{equation}

\begin{rmk}
A few comments about composition in a double category are in order.
A {\it compatible arrangement} is intuitively a pasting diagram of
squares in a double category. It was shown in \cite{dawsonpare1993}
that if a compatible arrangement has a composite, then this
composite does not depend on the order of composition, although
there may be compatible arrangements in a given double category that
do not admit a composite at all. We show in Corollary
\ref{compatiblearrangements} that all compatible arrangements in a
double category with folding admit a unique composite. We implicitly
use this existence and uniqueness throughout.
\end{rmk}

\begin{rmk}
The assignments $j \mapsto i^h_j$ and $f \mapsto i^v_f$ preserve
compositions. Preservation of units follows from the other axioms:
an application of the interchange law to the diagram of identity
morphisms
$$\xymatrix{A \ar[r] \ar[d] \ar@{}[dr]|{i^v_{1^h_A}} & A \ar[r] \ar[d]
\ar@{}[dr]|{i^h_{1^v_A}} & A \ar[d]  \\ A \ar[r] \ar[d]
\ar@{}[dr]|{i^h_{1^v_A}} & A \ar[r] \ar[d] \ar@{}[dr]|{i^v_{1^h_A}}
& A \ar[d]
\\ A \ar[r] & A \ar[r] & A }$$
shows that $i^h_{1^v_A}=i^v_{1^h_A}$. We abbreviate this identity
square with identity boundary simply by $i_A$. This proof does not
work for pseudo double categories, so $i^h_{1^v_A}=i^v_{1^h_A}$ is
an axiom in Definition \ref{pseudodoubledefn}.
\end{rmk}

\begin{rmk}
As a last comment about composition we remark that double categories
are not required to admit mixed compositions between horizontal and
vertical morphisms. Typically, horizontal and vertical morphisms are
different, as Example \ref{Rngs} shows.
\end{rmk}

\begin{defn} \label{horizontal2category}
Let $\mathbb{D}$ be a double category. Then $\mathbf{H}\mathbb{D}$
denotes the {\it horizontal 2-category} of $\mathbb{D}$. Its objects
are the objects of $\mathbb{D}$, its morphisms are the horizontal
morphisms of $\mathbb{D}$, and its 2-cells are the squares of
$\mathbb{D}$ which have vertical identities on the left and right
sides. The underlying 1-category of $\mathbf{H}\mathbb{D}$ is
denoted $(\mathbf{H}\mathbb{D})_0$. The {\it vertical 2-category}
$\mathbf{V}\mathbb{D}$ of $\mathbb{D}$ and its underlying 1-category
$(\mathbf{V}\mathbb{D})_0$ are defined analogously.
\end{defn}

\begin{defn}
If $\mathbb{D}$ and $\mathbb{E}$ are double categories, a {\it
double functor} $\xymatrix@1{F:\mathbb{D} \ar[r] & \mathbb{E}}$ is
an {\it internal functor in $Cat$}. This consists of functors
$\xymatrix@1{F_0:\mathbb{D}_0 \ar[r] & \mathbb{E}_0}$ and
$\xymatrix@1{F_1:\mathbb{D}_1 \ar[r] & \mathbb{E}_1}$  such that the
diagrams
$$\xymatrix@C=3pc@R=3pc{\mathbb{D}_1 \times_{\mathbb{D}_0}
\mathbb{D}_1 \ar[r] \ar[d]_{F_1
\times_{F_0} F_1} & \mathbb{D}_1 \ar[d]|{\tb{F_1}} & \mathbb{D}_0 \ar[l]_-{\eta} \ar[d]^{F_0} \\
\mathbb{E}_1 \times_{\mathbb{E}_0} \mathbb{E}_1 \ar[r]
& \mathbb{E}_1 & \mathbb{E}_0 \ar[l]^-{\eta} } $$
commute. In other words, a double functor consists of functions
$$\xymatrix{Obj \hspace{1mm} \mathbb{D} \ar[r] & Obj \hspace{1mm}  \mathbb{E}}$$
$$\xymatrix{Hor \hspace{1mm} \mathbb{D} \ar[r] & Hor \hspace{1mm} \mathbb{E}}$$
$$\xymatrix{Ver \hspace{1mm} \mathbb{D} \ar[r] & Ver \hspace{1mm} \mathbb{E}}$$
$$\xymatrix{Squares \hspace{1mm} \mathbb{D} \ar[r] & Squares \hspace{1mm} \mathbb{E}}$$
which preserve all sources, targets, compositions, and units.
\end{defn}

\begin{examp}
We can obtain double categories from a 2-category $\mathbf{C}$ in
several ways. The double category $\mathbb{H}\mathbf{C}$ has the
same objects as $\mathbf{C}$, horizontal morphisms are the morphisms
of $\mathbf{C}$, the vertical morphisms are all trivial, and the
squares are the 2-cells of $\mathbf{C}$. The double category
$\mathbb{V}\mathbf{C}$ is defined similarly, only this time all
horizontal morphisms are trivial. Any 2-functor
$\xymatrix@1{\mathbf{B} \ar[r] & \mathbf{C}}$ induces double
functors $\xymatrix@1{\mathbb{H}\mathbf{B} \ar[r] &
\mathbb{H}\mathbf{C}}$ and $\xymatrix@1{\mathbb{V}\mathbf{B} \ar[r]
& \mathbb{V}\mathbf{C}}$.
\end{examp}

\begin{examp} \label{quintetexample}
Another double category associated to a 2-category $\mathbf{C}$ is
Ehresmann's double category $\mathbb{Q}\mathbf{C}$ of {\it quintets
of $\mathbf{C}$}. Its objects are the objects of $\mathbf{C}$,
horizontal and vertical morphisms are the morphisms of $\mathbf{C}$,
and the squares $\alpha$ as in (\ref{sourcetarget}) are the 2-cells
$\xymatrix@1{\alpha:k \circ f \ar@{=>}[r] & g \circ j }$. Any
2-functor $\xymatrix@1{\mathbf{B} \ar[r] & \mathbf{C}}$ induces a
double functor $\xymatrix@1{\mathbb{Q}\mathbf{B} \ar[r] &
\mathbb{Q}\mathbf{C}}$. Note that the horizontal 2-category
$\mathbf{H}\mathbb{Q}\mathbf{C}$ is $\mathbf{C}$. The vertical
2-category $\mathbf{V}\mathbb{Q}\mathbf{C}$ is $\mathbf{C}$ with the
2-cells reversed. We could just as well have chosen our quintets to
consist of 2-cells $\xymatrix@1{g \circ j \ar@{=>}[r] &  k \circ f}$
instead, but then the roles of $\mathbf{H}\mathbb{Q}\mathbf{C}$ and
$\mathbf{V}\mathbb{Q}\mathbf{C}$ would be switched. In this article
we use the former convention because compatibility with $\mathbf{H}$
is important for folding. If $I$ is a 1-category viewed as a
2-category with only trivial 2-cells, then $\mathbb{Q}I$ is the {\it
double category $\Box I$ of commutative squares in I}. A boundary
admits a unique square if and only if the boundary is a commutative
square. The double categories $\mathbb{Q}\mathbf{C}$ are {\it
edge-symmetric} double categories as in \cite{brownmosa99} because
the horizontal and vertical edge categories coincide.
\end{examp}

\begin{examp} \label{adjunctions}
An {\it adjunction} in a 2-category $\mathbf{C}$ consists of two
morphisms $\xymatrix@1{j_1:A \ar[r] & C}$ and $\xymatrix@1{j_2:C
\ar[r] & A}$ and two 2-cells $\xymatrix@1{\eta:1_C \ar@{=>}[r] & j_1
\circ j_2}$ and $\xymatrix@1{\varepsilon:j_2\circ j_1 \ar@{=>}[r] &
1_A}$ which satisfy the familiar triangle identities. Here $j_1$ is
the {\it right} adjoint and  $j_2$ is the {\it left} adjoint. The
adjunctions in $\mathbf{C}$ form a double category
$\mathbb{A}\text{d}\mathbf{C}$ with objects the objects of
$\mathbf{C}$, horizontal morphisms the morphism of $\mathbf{C}$,
vertical morphisms the adjunctions in $\mathbf{C}$ (with direction
given by the right adjoint), and squares $\alpha$ as in
(\ref{sourcetarget}) the 2-cells $\xymatrix@1{\alpha:k_1 \circ f
\ar@{=>}[r] & g \circ j_1}$. This double category (with squares
reversed) was reviewed in \cite{kelly} to describe the sense in
which {\it mates under adjunctions} are compatible with composition
and identity. There is a forgetful double functor
$\xymatrix@1{\mathbb{A}\text{d}\mathbf{C} \ar[r] &
\mathbb{Q}\mathbf{C}}$. A related double category of certain {\it
adjoint squares} was introduced and studied in
\cite{palmquistthesis} and \cite{palmquist1971}.
\end{examp}

We will also have occasion to use double natural transformations.
There are two types: horizontal and vertical.

\begin{defn} \label{horizontalnaturaltransformation}
If $\xymatrix@1{F,G:\mathbb{D} \ar[r] & \mathbb{E}}$ are double
functors, then a {\it horizontal natural transformation}
$\xymatrix@1{\theta:F \ar@{=>}[r] & G}$ as in \cite{grandisdouble1}
assigns to each object $A$ a horizontal arrow $\xymatrix@1{\theta
A:FA \ar[r] & GA}$ and assigns to each vertical morphism $j$ a
square
$$\xymatrix@R=3pc@C=3pc{FA \ar[r]^{\theta A} \ar[d]_{Fj} \ar@{}[dr]|{\theta j} & GA \ar[d]^{Gj} \\
FC \ar[r]_{\theta C} & GC}$$ such that:
\begin{enumerate}
\item
For all $A \in \mathbb{D}$, we have $\theta 1^v_A =i^v_{\theta A }$,
\item
For composable vertical morphisms $j$ and $k$,
$$
\begin{array}{c}
\xymatrix@R=3pc@C=3pc{FA \ar[r]^{\theta A } \ar[d]_{F\vcomp{j}{k}}
\ar@{}[dr]|{\theta \vcomp{j}{k} } & GA \ar[d]^{F\vcomp{j}{k}} \\
FE \ar[r]_{\theta E} & GE}
\end{array}
=
\begin{array}{c}
\xymatrix@R=3pc@C=3pc{FA \ar[r]^{\theta A } \ar[d]_{Fj} \ar@{}[dr]|{\theta j } & GA \ar[d]^{Gj} \\
FC \ar[r]|{\lr{\theta C }} \ar[d]_{Fk} \ar@{}[dr]|{\theta k } & GC \ar[d]^{Gk} \\
FE \ar[r]_{\theta E} & GE,}
\end{array}
$$
\item
For all $\alpha$ as in (\ref{sourcetarget}),
$$
\begin{array}{c}
\xymatrix@R=3pc@C=3pc{FA \ar[r]^{Ff} \ar[d]_{Fj}
\ar@{}[dr]|{F\alpha} & FB \ar[r]^{\theta B} \ar[d]|{\tb{Fk}}
\ar@{}[dr]|{\theta k} & GB \ar[d]^{Gk} \\ FC \ar[r]_{Fg} & FC
\ar[r]_{\theta C} & GD}
\end{array}
=
\begin{array}{c}
\xymatrix@R=3pc@C=3pc{FA \ar[r]^{\theta A} \ar[d]_{Fj}
\ar@{}[dr]|{\theta j} & GA \ar[r]^{Gf} \ar[d]|{\tb{Gj}}
\ar@{}[dr]|{G\alpha} & GB \ar[d]^{Gk} \\ FC \ar[r]_{\theta C} & GC
\ar[r]_{Gg} & GD.}
\end{array}
$$
\end{enumerate}
\end{defn}
A horizontal natural transformation is the same as an {\it internal
natural transformation in $Cat$}. We also have vertical natural
transformations:

\begin{defn} \label{verticalnaturaltransformation}
If $\xymatrix@1{F,G:\mathbb{D} \ar[r] & \mathbb{E}}$ are double
functors, then a {\it vertical natural transformation}
$\xymatrix@1{\sigma:F \ar@{=>}[r] & G}$ as in \cite{grandisdouble1}
assigns to each object $A$ a vertical arrow $\xymatrix@1{\sigma A:FA
\ar[r] & GA}$ and assigns to each horizontal morphism $f$ a square
$$\xymatrix@R=3pc@C=3pc{FA \ar[d]_{\sigma A} \ar[r]^{Ff} \ar@{}[dr]|{\sigma f} & FB  \ar[d]^{\sigma B} \\
GA  \ar[r]_{Gf} & GB}$$ such that:
\begin{enumerate}
\item
For all objects $A \in \mathbb{D}$, we have $\sigma
1^h_A=i^h_{\sigma A}$,
\item
For all composable horizontal morphisms $f$ and $g$, $$\sigma [ f \
g  ]=[\sigma f \ \ \sigma g  ],$$
\item
For all $\alpha$ as in (\ref{sourcetarget}),
$$\vcomp{F \alpha}{\sigma g}=\vcomp{\sigma f}{G\alpha}.$$
\end{enumerate}
\end{defn}

\begin{examp}
An oplax natural transformation between 2-functors
$\xymatrix@1{\mathbf{B} \ar[r] & \mathbf{C}}$ is the same as a
vertical natural transformation between the induced double functors
$\xymatrix@1{\mathbb{H}\mathbf{B} \ar[r] & \mathbb{H}\mathbf{C}}$.
The components are necessarily trivial.
\end{examp}

To compare $I$-categories (and more generally pseudo $I$-categories)
with certain double categories, we extend Brown and Mosa's notion of
{\it folding} to non-edge-symmetric double categories. We prove that
a folding is equivalent to a {\it connection pair} in Lemmas
\ref{connectionpairtohorizontalization}-\ref{lastpiece} and Theorem
\ref{connection=horizontalization}. In the case of edge-symmetric
double categories, a connection pair (with trivial holonomy) is the
same as a {\it thin structure} as shown in \cite{brownmosa99}, and
in higher dimensions in \cite{higgins2005}.
Edge-symmetric foldings were used already in
\cite{brownhigginscubes} to prove that the category of crossed
complexes is equivalent to the category of cubical
$\omega$-groupoids, and were generalized to all dimensions
in \cite{alaglthesis}.
More recently, foldings found important
applications in \cite{alglbrownsteiner2002} and \cite{higgins2005}.
To define foldings, we recall Brown and Spencer's notion of {\it
holonomy} in \cite{brownspencer76}:

\begin{defn} \label{holonomy}
A {\it holonomy} for a double category $\mathbb{D}$ is a 2-functor
$\xymatrix@1{(\mathbf{V}\mathbb{D})_0 \ar[r] &
\mathbf{H}\mathbb{D}}$ which is the identity on objects. In other
words, a {\it holonomy} associates to a vertical morphism a
horizontal morphism with the same domain and range in a functorial
way.
\end{defn}

\begin{rmk} \label{holonomy=mixedcomposition}
If a double category is equipped with a holonomy, then we can define
a composition of vertical morphisms $j$ with morphisms and 2-cells
of $\mathbf{H}\mathbb{D}$ by
$$f \circ j := f \circ \overline{j}$$
$$j \circ g := \overline{j} \circ g$$
$$ \alpha \circ j :=\alpha \circ i^v_{\overline{j}}$$
$$j \circ \beta:= i^v_{\overline{j}} \circ \beta$$
to obtain morphisms and 2-cells of $\mathbf{H}\mathbb{D}$. Here
$\circ$ on the right-hand side designates horizontal composition of
morphisms and 2-cells of the 2-category $\mathbf{H}\mathbb{D}$.
These mixed compositions satisfy the obvious axioms: associativity,
unitality, and the usual axioms of left and right whiskering (see
\cite{streetstructures}). Conversely, a double category with a mixed
composition satisfying these axioms admits a holonomy defined by $j
\mapsto 1_C \circ j=j \circ 1_A$. These two procedures are inverse,
thus holonomies are the same as such mixed compositions.
\end{rmk}

\begin{rmk} \label{extendingD}
Given a double category $\mathbb{D}$ equipped with a holonomy, or
equivalently with mixed composition, one can construct a new double
category $\mathbb{D}'$ with an inclusion holonomy.  The objects and
vertical 1-categories of $\mathbb{D}$ and $\mathbb{D}'$ are the
same, while the set of horizontal morphisms of $\mathbb{D}'$ is the
disjoint union of the sets of horizontal and vertical morphisms of
$\mathbb{D}$. Composition of horizontal morphisms in $\mathbb{D}'$
is the mixed composition, with identities the included vertical
identities. The squares of $\mathbb{D'}$ are the squares of
$\mathbb{D}$ along with vertical identity squares for the horizontal
morphisms of $\mathbb{D}'$ which come from vertical morphisms of
$\mathbb{D}$. We equip $\mathbb{D'}$ with a holonomy by including
the vertical morphisms, so that the double functor
$\xymatrix@1{\mathbb{D}' \ar[r] & \mathbb{D}}$ preserves the
holonomies. We will apply this construction in Example
\ref{worldsheetexample} and Example \ref{ringsinY}.
\end{rmk}

\begin{defn} \label{foldingstructure}
A {\it folding on a double category $\mathbb{D}$} is a double
functor $\xymatrix@1{\Lambda:\mathbb{D} \ar[r] &
\mathbb{Q}\mathbf{H}\mathbb{D}}$ which is the identity on the
horizontal 2-category $\mathbf{H}\mathbb{D}$ of $\mathbb{D}$ and is
faithfully full on squares. More specifically, a folding consists of
a holonomy $\xymatrix@1{j \ar@{|->}[r] & \overline{j}}$ and
bijections $\Lambda^{f,k}_{j,g}$ from squares in $\mathbb{D}$ with
boundary
\begin{equation} \label{boundary1}
\xymatrix@R=3pc@C=3pc{A \ar[r]^f \ar[d]_j & B \ar[d]^k \\ C \ar[r]_g
& D}
\end{equation}
to squares in $\mathbb{D}$ with boundary
\begin{equation} \label{boundary2}
\xymatrix@R=3pc@C=3pc{A \ar[r]^{[f \ \overline{k}]} \ar[d]_{1_A^v} &
D \ar[d]^{1_D^v}
\\ A \ar[r]_{[\overline{j} \ g]} & D,}
\end{equation}
such that:
\begin{enumerate}
\item
$\Lambda$ is the identity if $j$ and $k$ are vertical identity
morphisms.
\item \label{Lambdahorizontal} \label{Lambdahorizontalequ}
$\Lambda$ preserves horizontal composition of squares, \ie
$$\Lambda\left(
\begin{array}{c}\xymatrix@R=4pc@C=4pc{A \ar[r]^{f_1} \ar[d]_j
\ar@{}[dr]|{\alpha} & B \ar[r]^{f_2} \ar[d]|{\tb{k}}
\ar@{}[dr]|{\beta} & C \ar[d]^{\ell} \\ D \ar[r]_{g_1} & E
\ar[r]_{g_2} & F}
\end{array}
 \right) \quad
=
\begin{array}{c}
\xymatrix@R=4pc@C=4pc{ A \ar[r]^{[f_1 \ f_2 \ \overline{\ell}]}
\ar[d]_{1^v_A} \ar@{}[dr]|{[i_{f_1}^v \ \Lambda(\beta)]} & F \ar[d]^{1^v_F}  \\
A \ar[r]|{[f_1 \ \overline{k} \ g_2]} \ar[d]_{1^v_A}
\ar@{}[dr]|{[\Lambda(\alpha) \ i_{g_2}^v]} & F \ar[d]^{1^v_F} \\ A
\ar[r]_{[\overline{j} \ g_1 \ g_2]} & F.}
\end{array}$$
\item \label{Lambdavertical} \label{Lambdaverticalequ}
$\Lambda$ preserves vertical composition of squares, \ie
$$\Lambda \left( \begin{array}{c} \xymatrix@R=4pc@C=4pc{A \ar[d]_{j_1} \ar[r]^f
\ar@{}[dr]|\alpha & B \ar[d]^{k_1} \\ C \ar[r]|{\lr{g}}
\ar@{}[dr]|\beta \ar[d]_{j_2} & D \ar[d]^{k_2} \\ E \ar[r]_h &
F,}\end{array} \right) \quad
 =
\begin{array}{c} \xymatrix@R=4pc@C=4pc{A \ar[r]^{[ f \
\overline{k}_1 \ \overline{k}_2]} \ar[d]_{1^v_A}
\ar@{}[dr]|{[\Lambda(\alpha) \ i_{\overline{k}_2}^v]} & F
\ar[d]^{1^v_E} \\
A \ar[r]|{[\overline{j}_1 \ g \ \overline{k}_2]}
\ar@{}[dr]|{[i_{\overline{j}_1}^v \ \Lambda(\beta)]} \ar[d]_{1^v_A}
& F \ar[d]^{1^v_E} \\ A \ar[r]_{[\overline{j}_1 \ \overline{j}_2 \
h] } & F.}
\end{array}
$$
\item \label{Lambdaidentity} \label{Lambdaidentityequ}
$\Lambda$ preserves identity squares, \ie
$$\Lambda\left(
\begin{array}{c}
\xymatrix@R=4pc@C=4pc{A \ar[r]^{1_A^h} \ar[d]_{j} \ar@{}[dr]|{i_j^h} & A \ar[d]^j \\
B \ar[r]_{1_B^h} & B}
\end{array} \right) \quad
=
\begin{array}{c}
\xymatrix@R=4pc@C=4pc{A \ar[r]^{[1^h_A \ \overline{j}]}
\ar[d]_{1^v_A} \ar@{}[dr]|{i^v_{\overline{j}}} & B \ar[d]^{1^v_B} \\
A \ar[r]_{[ \overline{j} \  1^h_B]}  & B.}
\end{array}
$$
\end{enumerate}
\end{defn}

\begin{defn} \label{morphismwithfolding}
Let $\mathbb{D}$ and $\mathbb{E}$ be double categories with folding.
A {\it morphism of double categories with folding}
$\xymatrix@1{F:\mathbb{D} \ar[r] & \mathbb{E}}$ is a double functor
such that
$$F(\overline{j})=\overline{F(j)}$$
$$F(\Lambda^{\mathbb{D}}(\alpha))=\Lambda^{\mathbb{E}}(F(\alpha))$$
for all vertical morphisms $j$ and squares $\alpha$ in $\mathbb{D}$.
This is a double functor $F$ such that
$$\xymatrix@C=3pc{\mathbb{D} \ar[r]^F \ar[d] & \mathbb{E} \ar[d] \\
\mathbb{Q}\mathbf{H}\mathbb{D} \ar[r]_{\mathbb{Q}\mathbf{H}F} &
\mathbb{Q}\mathbf{H}\mathbb{E}}$$ commutes.
\end{defn}

\begin{defn} \label{2cellwithfolding}
Let $\xymatrix@1{F,G:\mathbb{D} \ar[r] & \mathbb{E}}$ be morphisms
of double categories with folding.  A horizontal natural
transformation $\xymatrix@1{\theta:F \ar@{=>}[r] & G}$ is {\it
compatible with folding} if for all vertical morphisms $j$ the
following equation holds.
$$\Lambda\left(
\begin{array}{c}
\xymatrix@R=4pc@C=4pc{FA \ar[r]^{\theta A} \ar[d]_{Fj} \ar@{}[dr]|{\theta j} & GA \ar[d]^{Gj} \\
FC \ar[r]_{\theta C} & GC}
\end{array} \right) \quad
=
\begin{array}{c}
\xymatrix@R=4pc@C=4pc{FA \ar[r]^{[\theta A \ G\overline{j}]}
\ar[d]_{1^v_{FA}} \ar@{}[dr]|{i^v_{[\theta A \
G\overline{j}]}} & GC \ar[d]^{1^v_{GC}} \\
FA \ar[r]_{[ G\overline{j} \ \theta C]} \ar[r] & GC}
\end{array}
$$
A vertical natural transformation $\xymatrix@1{\sigma:F \ar@{=>}[r]
& G}$ is {\it compatible with folding} if for all vertical morphisms
$j$ the following equation holds.
$$\Lambda\left(
\begin{array}{c}
\xymatrix@R=4pc@C=4pc{FA \ar[r]^{F \overline{j}} \ar[d]_{\sigma A}
\ar@{}[dr]|{\sigma \overline{j}} & FC \ar[d]^{\sigma C} \\
GA \ar[r]_{G \overline{j}} & GC}
\end{array} \right) \quad
=
\begin{array}{c}
\xymatrix@R=4pc@C=4pc{FA \ar[r]^{[F \overline{j} \ \overline{\sigma
C}]} \ar[d]_{1^v_{FA}} \ar@{}[dr]|{i^v_{[F \overline{j} \ \sigma
C]}} & GC \ar[d]^{1^v_{GC}} \\
FA \ar[r]_{[\overline{\sigma A}  \ G\overline{j}]} \ar[r] & GC}
\end{array}
$$
\end{defn}

\begin{rmk} \label{2cellcompatiblewithfolding}
The compatibility of a horizontal natural transformation with
folding implies that it is entirely determined by its restriction to
the horizontal 2-category. Even more is true, any 2-natural
transformation between the underlying horizontal 2-functors of two
morphisms gives rise to a horizontal natural transformation
compatible with folding, since the compatibility defines $\theta j$
and the folding axioms guarantee (ii) and (iii) of Definition
\ref{horizontalnaturaltransformation}. The analogous remark for
vertical natural transformations does not hold, since compatibility
only concerns $\sigma \overline{j}$ and not the more general $\sigma
f$.
\end{rmk}

\begin{defn} \label{morphismoffoldingstructures}
Let $\mathbb{D}$ be a double category equipped with two foldings
$\xymatrix@1{\Lambda_1,\Lambda_2:\mathbb{D} \ar[r] &
\mathbb{Q}\mathbf{H}\mathbb{D}}.$ Then a {\it morphism of foldings}
$\xymatrix@1{\theta:\Lambda_1 \ar[r] & \Lambda_2}$ is a {\it
horizontal} natural transformation
$$\xymatrix{\theta:\Lambda_1 \ar@{=>}[r] & \Lambda_2}$$
with identity components. Equivalently, writing $j \mapsto
\overline{j}$ for the holonomy of $\Lambda_1$ and $j \mapsto
\overline{j}$ for the holonomy of $\Lambda_2$, a morphism of
foldings $\theta$ assigns to each vertical morphism $j$ a square
\begin{equation} \label{foldingmorphismdiagram}
\xymatrix@R=3pc@C=3pc{A \ar@{}[dr]|{\theta j} \ar[r]^{\overline{j}}
\ar[d]_{1_A^v} & C \ar[d]^{1_C^v}
\\ A \ar[r]_{\overline{\overline{j}}} & C }
\end{equation}
such that:
\begin{enumerate}
\item
$\theta$ preserves identities, \ie $\theta 1_A^v =i_A$.
\item
$\theta$ preserves compositions, \ie $\theta
\vcomp{j_1}{j_2}=[\theta j_1 \ \ \theta j_2 ]$.
\item
$\theta$ is natural, \ie
$$\begin{array}{c}
\xymatrix@R=3pc@C=3pc{A \ar[r]^{[f \ \overline{k}]} \ar[d]_{1^v_A}
\ar@{}[dr]|{\Lambda_1(\alpha)} & D \ar[d]^{1^v_D} \\ A
\ar[r]|{[\overline{j} \ g]} \ar[d]_{1^v_A} \ar@{}[dr]|{[\theta j \
i^v_g]} & D \ar[d]^{1^v_D} \\ A \ar[r]_{[\overline{\overline{j}} \
g]} & D}
\end{array}=
\begin{array}{c}
\xymatrix@R=3pc@C=3pc{A \ar[r]^{[f \ \overline{k}]} \ar[d]_{1^v_A}
\ar@{}[dr]|{[i^v_f \ \theta k ]} & D \ar[d]^{1^v_D} \\ A \ar[r]|{[f
\ \overline{\overline{k}}]} \ar[d]_{1^v_A}
\ar@{}[dr]|{\Lambda_2(\alpha)} & D \ar[d]^{1^v_D} \\ A
\ar[r]_{[\overline{\overline{j}} \  g]} & D.}
\end{array}$$
\end{enumerate}
\end{defn}
It may appear that a morphism of foldings is a vertical natural
transformation because of Diagram (\ref{foldingmorphismdiagram}).
But this is not so, since $\theta j$ is a square in the
$\mathbb{Q}\mathbf{H}\mathbb{D}$ with trivial horizontal components,
and such a square is precisely of the form
(\ref{foldingmorphismdiagram}). One could alternatively interpret
(\ref{foldingmorphismdiagram}) to be an oplax natural transformation
with identity components between the 2-functors that constitute the
holonomies, though the naturality of this comparison 2-cell is not
equivalent to the full naturality of (iii).

A double category with folding is determined by its vertical
1-catego-ry, its horizontal 2-category, and the holonomy.
Vice-a-versa, one can construct a double category with folding from
a 2-category equipped with a 2-functor resembling a holonomy. This
will be made precise in Theorem \ref{YZ}, which states the key
feature of foldings: the 2-category of double categories with
folding is 2-equivalent to the 2-category of certain 2-functors. The
pseudo counterpart of Theorem \ref{YZ} is Theorem \ref{YZpseudo}).

The squares of a double category with folding are determined by the
2-cells of the underlying horizontal 2-category via the folding.  A
folding {\it horizontalizes} a double category in the sense that it
maps a double category to its underlying horizontal 2-category in a
functorial way in terms of quintets. Thus, the quintessential
example of a double category with folding is the double category of
quintets of a 2-category as follows.

\begin{examp}
Let $\mathbf{C}$ be a 2-category and $\mathbb{Q}\mathbf{C}$ the
double category of quintets of $\mathbf{C}$ as in Example
\ref{quintetexample}. The holonomy is the inclusion of the vertical
1-category $\mathbf{C}_0$ into the horizontal 2-category
$\mathbf{C}$, and the folding maps are the identity: the squares of
$\mathbb{Q}\mathbf{C}$ with boundary (\ref{boundary1}) are by
definition the 2-cells in $\mathbf{C}$ with boundary
(\ref{boundary2}). In fact $\mathbb{Q}$ is a 2-functor from the
2-category of small 2-categories to the 2-category of double
categories with folding, morphisms, and horizontal natural
transformations compatible with folding. As a special case of
$\mathbb{Q}\mathbf{C}$, the double category $\Box I$ of commutative
squares in a 1-category $I$ admits a folding. The folding on $\Box
I$ is unique. In fact, we'll see in Theorem
\ref{foldingsareisomorphic} that foldings are unique up to
isomorphism.
\end{examp}

\begin{examp}
The double category $\mathbb{A}\text{d}\mathbf{C}$ of Example
\ref{adjunctions} admits a canonical folding: the holonomy sends an
adjunction to its right adjoint part.  The forgetful double functor
$\xymatrix@1{\mathbb{A}\text{d}\mathbf{C} \ar[r] &
\mathbb{Q}\mathbf{C}}$ is an example of a morphism of double
categories with folding.
\end{examp}

In Section \ref{pseudodoublesection} we will extend the notion of
folding to pseudo double categories. The extended notion has more
examples, such as the pseudo double categories $\mathbb{R}$ng and
$\mathbb{W}$ of bimodules and worldsheets.

Connection pairs on double categories can be found in
\cite{brownmosa99} and \cite{spencer77}. In the terminology of
\cite{grandisdouble2}, a connection pair is a functorial choice of a
so-called orthogonal companion for each vertical morphism.

\begin{defn} \label{connectionpair}
A {\it connection pair} on a double category consists of a holonomy
$\xymatrix@1{j \ar@{|->}[r] & \overline{j}}$ and an assignment of a
pair of squares
$$\xymatrix@R=3pc@C=3pc{A \ar[r]^{\overline{j}} \ar[d]_j \ar@{}[dr]|{\Gamma(j)}
& C \ar[d]^{1_C^v} & & A \ar[d]_{1_A^v} \ar[r]^{1_A^h}
 \ar@{}[dr]|{\Gamma'(j)} & A \ar[d]^j \\
C \ar[r]_{1_C^h} & C & & A \ar[r]_{\overline{j}} & C}$$ to each
vertical morphism $j$ such that:
\begin{enumerate}
\item \label{connectionpairidentity}
$\Gamma$ and $\Gamma'$ preserve identities.
$$\Gamma(1_A^v)=i_{A} \hspace{1.5in} \Gamma'(1_A^v)=i_{A}$$
\item \label{connectionpaircomposition}
$\Gamma$ and $\Gamma'$ preserve compositions, \ie the {\it transport laws} hold.
$$
\begin{array}{c}
\Gamma\left(\begin{bmatrix} j_1 \\ j_2
\end{bmatrix} \right)\quad =
\end{array}
\begin{array}{l}
\xymatrix@R=3pc@C=3pc{\ar[r]^{\overline{j}_1} \ar[d]_{j_1}
\ar@{}[dr]|{\Gamma(j_1)} & \ar[r]^{\overline{j}_2}
 \ar[d] \ar@{}[dr]|{i_{\overline{j}_2}^v} & \ar[d]  \\
\ar[r] \ar[d]_{j_2} \ar@{}[dr]|{i_{j_2}^h} &
\ar[r]|{\lr{\overline{j}_2}} \ar[d]|{\tb{j_2}}
\ar@{}[dr]|{\Gamma(j_2)} & \ar[d]
\\\ar[r]_{\phantom{j_1}} & \ar[r]  &  }
\end{array}
$$
$$
\begin{array}{c}
\Gamma'\left(\begin{bmatrix} j_1 \\ j_2
\end{bmatrix} \right)\thickspace\thinspace =
\end{array}
\begin{array}{l}
\xymatrix@R=3pc@C=3pc{\ar[r]^{\phantom{j_1}} \ar[d]
\ar@{}[dr]|{\Gamma'(j_1)} & \ar[r]
 \ar[d]|{\tb{j_1}} \ar@{}[dr]|{i_{j_1}^h} & \ar[d]^{j_1}  \\
\ar[r]|{\lr{\overline{j}_1}} \ar[d]
\ar@{}[dr]|{i_{\overline{j}_1}^v} & \ar[r] \ar[d]
\ar@{}[dr]|{\Gamma'(j_2)} & \ar[d]^{j_2}
\\ \ar[r]_{\overline{j}_1} & \ar[r]_{\overline{j}_2}  &  }
\end{array}
$$
(unlabelled arrows are the identities).
\item \label{three}
$$
i^v_{\overline{j}}\thickspace\thinspace= \begin{array}{c}
\xymatrix@R=3pc@C=3pc{ \ar[r] \ar[d] \ar@{}[dr]|{\Gamma'(j)} &
\ar[r]^{\overline{j}} \ar[d]|{\tb{j}} \ar@{}[dr]|{\Gamma(j)}
 & \ar[d] \\  \ar[r]_{\overline{j}} & \ar[r] & }\end{array} \hspace{.5in}
 i^h_j\thickspace\thickspace\thickspace=\begin{array}{c}
\xymatrix@R=3pc@C=3pc{\ar[d] \ar[r] \ar@{}[dr]|{\Gamma'(j)} &  \ar[d]^{j} \\
\ar[r]|{\lr{\overline{j}}} \ar@{}[dr]|{\Gamma(j)} \ar[d]_{j} &  \ar[d] \\
 \ar[r] & } \end{array}$$
(unlabelled arrows are the identities).
\end{enumerate}
\end{defn}

We now work towards a proof of Theorem
\ref{connection=horizontalization}, which states that the data for a
connection pair is equivalent to the data for a folding. This proof
is essentially a slight generalization of an argument in Section 5
of \cite{brownmosa99}. The idea goes back to Spencer in
\cite{spencer77} and to the quintets of Ehresmann in
\cite{ehresmannquintett63}. If a double category admits a folding,
then that folding is unique up to isomorphism.






\begin{lem} \label{connectionpairtohorizontalization}
If $(\Gamma,\Gamma')$ is a connection pair on a double category,
then $$\Lambda^{f,k}_{j,g}(\alpha):=\begin{bmatrix} \Gamma'(j) &
\alpha & \Gamma(k)
\end{bmatrix}$$ defines a folding.
\end{lem}
\begin{pf}
First we show that the holonomy is part of a double functor
$\xymatrix@1{D \ar[r] & \mathbb{Q}\mathbf{H}\mathbb{D}}$.
\begin{enumerate}
\item
Follows from Definition \ref{connectionpair}
(\ref{connectionpairidentity}).
\item \label{compositionverification}
In the notation of Definition \ref{connectionpair}, we have
$$\aligned \Lambda(\begin{bmatrix} \alpha & \beta \end{bmatrix})&= \begin{bmatrix} \Gamma'(j) &
\alpha & \beta & \Gamma(\ell)
\end{bmatrix} \\
&= \begin{bmatrix} i_1^v & i_{f_1}^v & \Gamma'(k) & \beta &
\Gamma(\ell)
\\ \Gamma'(j) & \alpha & \Gamma(k) & i^v_{g_2} & i_1^v \end{bmatrix} \\
&= \begin{bmatrix} i^v_{f_1} & \Lambda(\beta) \\
\Lambda(\alpha) & i^v_{g_2} \end{bmatrix}
\endaligned$$
by way of Definition \ref{connectionpair} (\ref{three}).
\item
This proof is similar to the proof just given in
(\ref{compositionverification}).
\item
$$\aligned \Lambda(i^h_j) &= \begin{bmatrix} \Gamma'(j) & i^h_j & \Gamma(j) \end{bmatrix} \\
&= \begin{bmatrix} \Gamma'(j) &  \Gamma(j)
\end{bmatrix} \\
&=i^v_{\overline{j}} \endaligned$$ by Definition
\ref{connectionpair} (\ref{three}).
\end{enumerate}

The double functor $\Lambda$ is surjective on squares, since
$$\xymatrix@R=3pc@C=3pc{\ar[r]^f \ar[d] \ar@{}[dr]|{i^v_f} & \ar[r] \ar[d] \ar@{}[dr]|{\Gamma'(k)}  & \ar[d]^k
\\ \ar[r]|{\lr{f}} \ar[d] & \ar[r]|{\lr{\overline{k}}} \ar@{}[d]|{\delta} &
\ar[d] \\ \ar[d]_j \ar[r]|{\lr{\overline{j}}} \ar@{}[dr]|{\Gamma(j)}
& \ar[r]|{\lr{g}}
\ar[d] \ar@{}[dr]|{i^v_g} & \ar[d] \\
\ar[r] & \ar[r]_g & }$$ maps to $\delta$ by the left half of
Definition \ref{connectionpair} (iii). The right half of Definition
\ref{connectionpair} (iii) shows that
$$\xymatrix@R=3pc@C=3pc{\ar[rr]^f \ar[d] & \ar@{}[d]|{i^v_f} & \ar[r] \ar[d] \ar@{}[dr]|{\Gamma'(k)} & \ar[d]^k
\\ \ar[d] \ar[r] \ar@{}[dr]|{\Gamma'(j)} & \ar[d]|{\tb{j}} \ar[r]|{\lr{f}} \ar@{}[dr]|{\alpha} &
\ar[d]|{\tb{k}} \ar[r]|{\lr{\overline{k}}} \ar@{}[dr]|{\Gamma(k)} &
\ar[d]
\\ \ar[d]_{j} \ar[r]|{\lr{\overline{j}}} \ar@{}[dr]|{\Gamma(j)} & \ar[d] \ar[r]|{\lr{g}}
& \ar[r] \ar@{}[d]|{i^v_g} & \ar[d]
\\ \ar[r] & \ar[rr]_g & &   }$$
equals $\alpha$, so that $\Lambda$ is injective on squares.
\end{pf}

We now prove the converse to Lemma
\ref{connectionpairtohorizontalization}:

\begin{lem} \label{horizontalizationtoconnectionpair}
If $\Lambda$ is a folding on a double category, then
$$\Gamma(j):=(\Lambda^{\overline{j},1}_{j,1})^{-1}(i^v_{\overline{j}}) \hspace{.5in}
\Gamma'(j):=(\Lambda^{1,j}_{1,\overline{j}})^{-1}(i^v_{\overline{j}})$$
defines a connection pair.
\end{lem}
\begin{pf}
\begin{enumerate}
\item
Follows because $\Lambda$ is the identity on squares with identity
boundary.
\item
An application of $\Lambda^{\overline{[j_1
j_2]},1}_{\vcomp{j_1}{j_2},1}$ to the right side of the equation for
$\Gamma$ in (\ref{connectionpaircomposition}) of Definition
\ref{connectionpair} yields the following.
$$\aligned \Lambda^{\overline{[j_1
j_2]},1}_{\vcomp{j_1}{j_2},1}\left(
\begin{bmatrix}
\Gamma(j_1) & i^v_{\overline{j}_2}
\\ i^h_{j_2} & \Gamma(j_2)
\end{bmatrix} \right) &=
\Lambda\left( \begin{bmatrix} \begin{bmatrix} \Gamma(j_1) & i^v_{\overline{j}_2} \end{bmatrix} \vspace{1mm} \\
\begin{bmatrix} i^h_{j_2} & \Gamma(j_2) \end{bmatrix}
\end{bmatrix} \right) \\
&= \begin{array}{c} \xymatrix@C=6pc@R=4pc{\ar[r]|{[\overline{j}_1 \
\overline{j}_2]} \ar[d] \ar@{}[dr]|{\Lambda([\Gamma(j_1) \
i^v_{\overline{j}_2}]) } & \ar[d] \\ \ar[d] \ar[r]|{[\overline{j}_1
\ \overline{j}_2]} \ar@{}[dr]|{i^v_{\overline{j}_1}
\Lambda([i^h_{j_2} \ \Gamma(j_2) ])} & \ar[d] \\
\ar[r]|{[\overline{j}_1 \ \overline{j}_2]} & } \end{array} \\
&= \left[ \begin{matrix}
\begin{bmatrix} i_{\overline{j}_1}^v & \Lambda(i^v_{\overline{j}_2}) \end{bmatrix} \vspace{1mm} \\
\begin{bmatrix} \Lambda(\Gamma(j_1)) & i_{\overline{j}_2}^v \end{bmatrix} \vspace{1mm} \\
\begin{bmatrix} i^v_{\overline{j}_1} & \Lambda(\Gamma(j_2)) \end{bmatrix} \vspace{1mm} \\
\begin{bmatrix} i^v_{\overline{j}_1} & \Lambda(i^h_{j_2}) \end{bmatrix}  \\
\end{matrix} \right] \\
&=\begin{matrix} \phantom{i^v_{[\overline{j}_1 \overline{j}_2]}} &
i^v_{[\overline{j}_1 \ \overline{j}_2]} &
\phantom{i^v_{[\overline{j}_1  \overline{j}_2]}}
\end{matrix}
\endaligned$$
From this we conclude $\Gamma \left(\begin{bmatrix} j_1 \\ j_2
\end{bmatrix} \right)
=
\begin{bmatrix}
\Gamma(j_1) & i^v_{\overline{j}_2}
\\ i^h_{j_2} & \Gamma(j_2)
\end{bmatrix}.$
A similar argument works for $\Gamma'$.
\item
An application of $\Lambda^{\overline{j},1}_{1,\overline{j}}$ to
$[\Gamma'(j) \hspace{2mm} \Gamma(j)]$ in (\ref{three}) of Definition
\ref{connectionpair} yields the following.
$$\aligned
\begin{bmatrix}
\Gamma'(j) & \Gamma(j) \end{bmatrix} &=
\Lambda^{\overline{j},1}_{1,\overline{j}} \left( \begin{bmatrix}
\Gamma'(j) & \Gamma(j) \end{bmatrix} \right) \\
&= \begin{bmatrix} \Lambda(\Gamma(j)) \\
\Lambda(\Gamma'(j)) \end{bmatrix} \\
&= \begin{bmatrix} i^v_{\overline{j}} \\
i^v_{\overline{j}} \end{bmatrix} \\
&= i^v_{\overline{j}}
\endaligned$$
A similar argument shows $\begin{bmatrix} \Gamma'(j)
\\ \Gamma(j) \end{bmatrix}=i^h_j$.
\end{enumerate}
\end{pf}

\begin{lem}
Let $(\Gamma,\Gamma')$ be a connection pair on a double category
with associated folding $\Lambda$ as in Lemma
\ref{connectionpairtohorizontalization}. Then the connection pair
associated to $\Lambda$ as in Lemma
\ref{horizontalizationtoconnectionpair} is the connection pair we
started with.
\end{lem}
\begin{pf}
By Definition \ref{connectionpair} (\ref{three}) we see that
$$(\Lambda^{\overline{j},1}_{j,1})(\Gamma(j))=\begin{bmatrix} \Gamma'(j)
& \Gamma(j) & i_1^v \end{bmatrix} =i^v_{\overline{j}}$$
$$(\Lambda^{1,j}_{1,\overline{j}})(\Gamma'(j))=\begin{bmatrix} i_1^v
& \Gamma'(j) & \Gamma(j) \end{bmatrix} =i^v_{\overline{j}}.$$
\end{pf}

\begin{lem} \label{lastpiece}
Let $\Lambda$ be a folding on a double category with associated
connection pair $(\Gamma,\Gamma')$ as in Lemma
\ref{horizontalizationtoconnectionpair}. Then the folding associated
to $(\Gamma,\Gamma')$ as in Lemma
\ref{connectionpairtohorizontalization} is the folding we started
with.
\end{lem}
\begin{pf}
The square $\begin{bmatrix} \Gamma'(j) & \alpha & \Gamma(k)
\end{bmatrix}$ has trivial vertical edges, and is therefore
preserved by $\Lambda$ as in (i) of Definition
\ref{foldingstructure}.
$$\aligned \begin{bmatrix} \Gamma'(j) & \alpha & \Gamma(k) \end{bmatrix}
&= \Lambda \left ( \begin{bmatrix} \Gamma'(j) & \alpha & \Gamma(k)
\end{bmatrix}\right ) \\
&= \Lambda\left( \begin{bmatrix} \begin{bmatrix}
\Gamma'(j) &\alpha \end{bmatrix} & \Gamma(k) \end{bmatrix}  \right) \\
&=\begin{bmatrix}
\begin{bmatrix}
i^v_f & \Lambda(\Gamma(k))
\end{bmatrix} \vspace{1mm} \\
\begin{bmatrix}
\Lambda(\begin{bmatrix} \Gamma'(j) & \alpha  \end{bmatrix})
\end{bmatrix}
\end{bmatrix} \\
&= \begin{bmatrix}
\begin{bmatrix}
i^v_f & \phantom{\Lambda(\alpha)}  & i^v_{\overline{k}}
\end{bmatrix}\vspace{1mm} \\
\begin{bmatrix}
\Lambda(\alpha)
\end{bmatrix}\vspace{1mm}\\
\begin{bmatrix}
\Lambda(\Gamma'(j)) & i_g
\end{bmatrix}\vspace{1mm}\\
\end{bmatrix} \\
&=\begin{bmatrix}
\begin{bmatrix}
i^v_f & \phantom{\Lambda(\alpha)}  & i^v_{\overline{k}}
\end{bmatrix}\vspace{1mm} \\
\begin{bmatrix}
\Lambda(\alpha)
\end{bmatrix}\vspace{1mm}\\
\begin{bmatrix}
i^v_{\overline{j}}&\phantom{\Lambda(\alpha)} & i_g
\end{bmatrix}\vspace{1mm}\\
\end{bmatrix} \\
&=\phantom{[[\Lambda(\alpha)}\Lambda(\alpha). \phantom{\Lambda(\alpha)}
\endaligned
$$
\end{pf}

\begin{thm} \label{connection=horizontalization}
The notions of connection pair and folding on a double category are
equivalent.
\end{thm}
\begin{pf}
This follows from Lemmas
\ref{connectionpairtohorizontalization}-\ref{lastpiece}.
\end{pf}

\begin{cor} \label{compatiblearrangements}
Any compatible arrangement in a double category with folding admits
a unique composite.
\end{cor}
\begin{pf}
We imitate the proof of the edge-symmetric case in
\cite{brownmosa99}. Theorem 4.1 and Theorem 5.1 in
\cite{dawsonpare1993} provide a useful criterion for every
compatible arrangement of a double category $\mathbb{D}$ to admit a
unique composite.  Suppose that every square $\alpha$ as in
(\ref{sourcetarget}) satisfies the following condition. If either
the horizontal source $j$ or the horizontal target $k$ admits a
(vertical) factorization, then that factorization extends to a
vertical factorization
$$\alpha=\vcomp{\alpha_1}{\alpha_2}.$$
In this situation, every compatible arrangement of $\mathbb{D}$
admits a unique composite.

We claim that any double category $\mathbb{D}$ with folding
satisfies this criterion. Let $(\Gamma,\Gamma')$ be the connection
pair associated to the folding. If $j=\vcomp{j_1}{j_2}$, then
$$\alpha \quad = \begin{array}{c}
\xymatrix@R=3pc@C=3pc{\ar@{=}[r] \ar[d]_{j_1} \ar@{}[dr]|{i^h_{j_1}}
& \ar[r]^f \ar@{}[ddr]|\alpha \ar[d]|{\tb{j_1}} & \ar[dd]^k \\
\ar@{=}[d] \ar@{=}[r] \ar@{}[dr]|{\Gamma'(j_2)} & \ar[d]|{\tb{j_2}} & \\
\ar[r]|{\lr{\overline{j_2}}} \ar[d]_{j_2} \ar@{}[dr]|{\Gamma(j_2)}
& \ar[r]|{\lr{g}} \ar@{=}[d] \ar@{}[dr]|{i^v_g} & \ar@{=}[d] \\
\ar@{=}[r] & \ar[r]_g & } \end{array}$$ is a vertical factorization
of $\alpha$ extending the factorization of $j$. A similar proof
works for factorizations of $k$.
\end{pf}

If a double category admits a folding, then the folding is
essentially unique:

\begin{thm} \label{foldingsareisomorphic}
Any two foldings on a double category $\mathbb{D}$ are isomorphic.
\end{thm}
\begin{pf}
Suppose $(\Lambda_1, j \mapsto \overline{j})$ and $(\Lambda_2, j
\mapsto \overline{\overline{j}})$ are foldings on $\mathbb{D}$ with
respective associated connection pairs $(\Gamma_1,\Gamma'_1)$ and
$(\Gamma_2,\Gamma'_2)$. We define a morphism
$\xymatrix@1{\theta:(\Lambda_1, j \mapsto \overline{j}) \ar[r] &
(\Lambda_2, j \mapsto \overline{\overline{j}})}$ of foldings by
$\theta j:=\Lambda_2(\Gamma_1(j))$. This is natural for $\alpha$ as
in (\ref{sourcetarget}) because
$$\aligned \begin{bmatrix} \Lambda_1(\alpha) \\
\begin{bmatrix} \Lambda_2(\Gamma_1(j)) & i^v_g\end{bmatrix} \end{bmatrix}
&= \begin{bmatrix} i^v_1 & \Gamma'_1(j) & \alpha & \Gamma_1(k)
\\ \Gamma'_2(j) & \Gamma_1(j) & i^v_g & i^v_1 \end{bmatrix} \\
&= \begin{bmatrix} \Gamma'_2(j) & \alpha & \Gamma_1(k) \end{bmatrix} \\
&= \begin{bmatrix} i^v_1 & i^v_f & \Gamma'_2(k) &  \Gamma_1(k)
\\ \Gamma'_2(j) & \alpha & \Gamma_2(k) & i^v_1 \end{bmatrix} \\
&= \begin{bmatrix} \begin{bmatrix} i^v_f & \Lambda_2(\Gamma_1(k))
\end{bmatrix} \\ \Lambda_2(\alpha) \end{bmatrix}.
\endaligned$$

An inverse to $\theta j$ is given by $\theta^{-1} j
:=\Lambda_2(\Gamma'_1(j))$.

$$\aligned \begin{array}{c}
\xymatrix@R=4pc@C=4pc{\ar[r]^{\overline{j}} \ar[d] \ar@{}[dr]|{\theta j} & \ar[d] \\
\ar[r]|{\lr{\overline{\overline{j}}}} \ar[d] \ar@{}[dr]|{\theta^{-1}
j} & \ar[d] \\ \ar[r]_{\overline{j}} & }\end{array} &=
\begin{array}{c}
\xymatrix@R=4pc@C=4pc{\ar[r]^{\overline{j}} \ar[d] \ar@{}[dr]|{\Lambda_2(\Gamma_1(j))} & \ar[d] \\
\ar[r]|{\lr{\overline{\overline{j}}}} \ar[d]
\ar@{}[dr]|{\Lambda_2(\Gamma'_1(j))} & \ar[d] \\
\ar[r]_{\overline{j}} & }
\end{array} \\ &=\Lambda_2(\begin{bmatrix}\Gamma'_1(j) & \Gamma_1(j)\end{bmatrix})\\
&=\Lambda_2(i^v_{\overline{j}}) \\
&=i^v_{\overline{j}} \\
\endaligned$$

The other direction $\vcomp{\theta^{-1}j}{\theta
j}=i^v_{\overline{\overline{j}}}$ follows similarly from
$$\begin{bmatrix} \Gamma'_1 (j) \\ \Gamma_1(j) \end{bmatrix}=i^h_{j}.$$
\end{pf}

A related structure on an edge-symmetric double category is a thin
structure as in \cite{brownmosa99}:

\begin{defn}
Let $\mathbb{D}$ be an edge-symmetric double category. Then a {\it
thin structure} on $\mathbb{D}$ is a double functor
$$\xymatrix{\Theta:\Box (\mathbf{H}\mathbb{D})_0 \ar[r] & \mathbb{D}}$$
which is the identity on objects and morphisms. Here $\Box
(\mathbf{H}\mathbb{D})_0$ is the double category of commutative
squares of morphisms of $\mathbb{D}$. The squares of $\mathbb{D}$ in
the image of $\Theta$ are called {\it thin}. Clearly, any
commutative boundary in $\mathbb{D}$ has a unique thin filler and
any composition of thin squares is thin.
\end{defn}

\begin{thm}
(Brown-Mosa in \cite{brownmosa99}) A thin structure and a connection
pair with trivial holonomy on an edge-symmetric double category
determine each other.
\end{thm}

\begin{cor} \label{folding=thin}
The notions of folding with trivial holonomy and thin structure on
an edge-symmetric double category are equivalent.
\end{cor}

After introducing double categories, foldings, connection pairs, and
thin structures in this section, we put them to use in an alternate
description of $I$-categories in the next section.

\section{$I$-Categories and Double Categories with Folding}

Although an $I$-category is not the same thing as an internal
category in $Cat$, it is of course a related concept. In this
section we show how $I$-categories are related to double categories
with folding in an explicit and elementary way. Surprisingly, both
notions are equivalent to the simple notion of a strict 2-functor
$\xymatrix@1{I \ar[r] & \mathcal{C}}$ that is the identity on
objects. We introduce three 2-categories $\mathcal{X}^{strict},
\mathcal{Y}^{strict}, \mathcal{Z}^{strict}$ and show that they are
2-equivalent if $I$ is a groupoid.  We also prove that
$\mathcal{Y}^{strict}$ and $\mathcal{Z}^{strict}$ remain
2-equivalent even if $I$ is merely a category. Unless stated
otherwise, $I$ denotes a fixed category.

\begin{notation} \label{Xnotation}
Let $\mathcal{X}^{strict}$ denote the 2-category of $I$-categories
as defined in Section \ref{strictcategories}. The morphisms in
$\mathcal{X}^{strict}$ are strict.
\end{notation}

\begin{notation} \label{Ynotation}
Let $\mathcal{Y}^{strict}$ denote the 2-category whose objects are
double categories $\mathbb{D}$ with folding such that
$(\mathbf{V}\mathbb{D})_0=I$. A morphism in $\mathcal{Y}^{strict}$
is a morphism of double categories with folding which is also the
identity on $(\mathbf{V}\mathbb{D})_0$.

The 2-cells of $\mathcal{Y}^{strict}$ are vertical natural
transformations that are compatible with folding and also have
identity components. More precisely, recall that
$\mathbf{H}\mathbb{D}(A,B)$ denotes the category whose objects are
the horizontal morphisms from $A$ to $B$ in $\mathbb{D}$ and whose
morphisms are 2-cells of $\mathbf{H}\mathbb{D}$ with source and
target such horizontal morphisms. If $\xymatrix@1{F:\mathbb{D}
\ar[r] & \mathbb{E}}$ is a morphism, we denote its restriction to
$\mathbf{H}\mathbb{D}(A,B)$ by
$\xymatrix@1{\mathbf{H}F_{A,B}:\mathbf{H}\mathbb{D}(A,B) \ar[r] &
\mathbf{H}\mathbb{E}(FA,FB)=\mathbf{H}\mathbb{E}(A,B)}$. If
$\xymatrix@1{F,G:\mathbb{D} \ar[r] & \mathbb{E}}$ are morphisms in
$\mathcal{Y}^{strict}$, then a 2-cell $\xymatrix@1{\sigma:F
\ar@{=>}[r] & G}$ in $\mathcal{Y}^{strict}$ assigns to each pair
$(A,B) \in I^2$ a natural transformation
$\xymatrix@1{\sigma_{A,B}:\mathbf{H}F_{A,B} \ar@{=>}[r] &
\mathbf{H}G_{A,B}}$ such that
$$\sigma^{\overline{j}}_{A,C}=i_{\overline{j}}^v$$
$$[ \ \sigma^f_{A,B} \ \ \sigma^g_{B,C}\ ]= \sigma^{[f \ g]}_{A,C} $$
$$\sigma_{A,A}^{1^h_A}=i_{1_A^h}^v$$
for all vertical morphisms $\xymatrix@1{j:A \ar[r] & C}$, composable
horizontal morphisms $f$,$g$, and all objects $A$. With these
definitions, $\mathcal{Y}^{strict}$ is a 2-category.
\end{notation}

\begin{notation} \label{Znotation}
Let $\mathcal{Z}^{strict}$ denote the 2-category of 2-categories
$\mathcal{C}$ with object set $Obj \hspace{1mm} I$ and equipped with
a strict 2-functor $\xymatrix@1{P:I \ar[r] & \mathcal{C}}$ from the
fixed category $I$ to $\mathcal{C}$ which is the identity on
objects.

A morphism from $\xymatrix@1{P:I \ar[r] & \mathcal{C}}$ to
$\xymatrix@1{P':I \ar[r] & \mathcal{C}'}$ in $\mathcal{Z}^{strict}$
is a strict 2-functor $\xymatrix@1{F:\mathcal{C} \ar[r] &
\mathcal{C}'}$ such that
$$\xymatrix@C=4pc@R=3pc{I \ar[r]^P \ar[dr]_{P'} & \mathcal{C}
\ar[d]^F \\ & \mathcal{C'}}$$ strictly commutes. We see that any
morphism $F$ is the identity on objects.

If $\xymatrix@1{F,G:P \ar[r] & P'}$ are morphisms in
$\mathcal{Z}^{strict}$, then a 2-cell $\xymatrix@1{\sigma:F
\ar@{=>}[r] & G}$ is {\it not} a 2-natural transformation from $F$
to $G$, but instead consists of natural transformations
$\xymatrix@1{\sigma_{A,B}:F_{A,B} \ar@{=>}[r] & G_{A,B}}$ for all
$A,B \in Obj \hspace{1mm} \mathcal{C}$ such that
$$\sigma_{A,C}^{P(j)}=i_{P'(j)}$$
$$\sigma^g_{B,C} \circ \sigma^f_{A,B}=\sigma^{g \circ f}_{A,C}$$
$$\sigma_{A,A}^{1^h_A}=i_{1_A^h}$$
for all $j \in I(A,C), f \in \mathcal{C}(A,B), g \in
\mathcal{C}(B,C),$ and objects $A$. Here $i_{P'(j)}$ denotes the
identity 2-cell on the morphism $P'(j)$ in the 2-category
$\mathcal{C}'$. The notation $\circ$ denotes the horizontal
composition of 2-cells as well as the composition of morphisms.

The 2-category $\mathcal{Z}^{strict}$ is similar to the 2-category
of pseudo 2-algebras over a theory in \cite{fiore3}.
\end{notation}

\begin{rmk}
If $I$ is a discrete category, then the objects and morphisms of
$\mathcal{X}^{strict}$, $\mathcal{Y}^{strict}$, and
$\mathcal{Z}^{strict}$ are simply 2-categories with object set $I$
and 2-functors that are the identity on objects. The 2-cells are
oplax natural transformations with identity components, which are
better viewed as vertical natural transformations with identity
components as in Section \ref{doublesection}.  We extend this
identification to general groupoids $I$ in Theorem \ref{XY} and
Theorem \ref{XZ} below, and then also to the weak situation in
Theorem \ref{XYpseudo} and Theorem \ref{XZpseudo}.
\end{rmk}

\begin{defn}
A 2-functor $\xymatrix@1{F:\mathcal{C} \ar[r] & \mathcal{D}}$ is a
2-equivalence if there exists a 2-functor $\xymatrix@1{G:\mathcal{D}
\ar[r] & \mathcal{C}}$ and 2-natural isomorphisms $1_\mathcal{C}
\cong GF$ and $FG \cong 1_{\mathcal{D}}$. The notion of
2-equivalence is the same as equivalence of $Cat$-enriched
categories, \ie a 2-functor which is surjective on objects up to
isomorphism and locally an isomorphism of categories.
\end{defn}

First we compare $\mathcal{Y}^{strict}$ and $\mathcal{Z}^{strict}$.
This 2-equivalence is the essential feature of foldings:

\begin{thm} \label{YZ}
Let $I$ be a category. The 2-category $\mathcal{Y}^{strict}$ of
double categories $\mathbb{D}$ with folding such that
$(\mathbf{V}\mathbb{D})_0=I$ is 2-equivalent to the 2-category
$\mathcal{Z}^{strict}$ of strict 2-functors $\xymatrix@1{P:I \ar[r]
& \mathcal{C}}$ that are the identity on objects as in Notation
\ref{Ynotation} and Notation \ref{Znotation}.
\end{thm}
\begin{pf}
We define two strict 2-functors
$\xymatrix@1{\mathcal{L}:\mathcal{Y}^{strict} \ar[r] &
\mathcal{Z}^{strict}}$ and
$\xymatrix@1{\mathcal{M}:\mathcal{Z}^{strict} \ar[r] &
\mathcal{Y}^{strict}}$ and show that $\mathcal{M}$ is surjective on
objects up to isomorphism and locally an isomorphism.

From $\mathbb{D} \in \mathcal{Y}^{strict}$, we define
$\mathcal{C}:=\mathbf{H}\mathbb{D}$ and take the functor
$\xymatrix@1{P:I \ar[r] & \mathcal{C}}$ to be the holonomy, in other
words $\mathcal{L}(\mathbb{D}):=P$. For a morphism $F$ and a 2-cell
$\sigma$ in $\mathcal{Y}^{strict}$, $\mathcal{L}(F):=\mathbf{H}(F)$
and $\mathcal{L}\sigma:=\sigma$.

For a strict 2-functor $\xymatrix@1{P:I \ar[r] & \mathcal{C}}$ in
$\mathcal{Z}^{strict}$, the double category $\mathcal{M}(P)$ has
vertical 1-category $I$ and horizontal 2-category $\mathcal{C}$. The
squares $\alpha$ of the form (\ref{sourcetarget}) are 2-cells
$\xymatrix@1{P(k) \circ f \ar@{=>}[r] & g \circ P(j)}$ in
$\mathcal{C}$. The holonomy is $P$ and the folding bijection is the
identity. Horizontal and vertical composition of squares, and the
respective units, are defined by the folding axioms in Definition
\ref{foldingstructure}. For a morphism $F$ in
$\mathcal{Z}^{strict}$, the double functor $\mathcal{M}(F)$ is the
identity on the vertical 1-category and $F$ on the horizontal
2-category. This in fact also defines $F$ on squares. Lastly
$\mathcal{M}(\sigma)=\sigma$ for all 2-cells $\sigma$ in
$\mathcal{Z}^{strict}$.

The 2-functor $\mathcal{M}$ is surjective on objects up to
isomorphism, since $\mathcal{M}\mathcal{L}\mathbb{D} \cong
\mathbb{D}$ for all $\mathbb{D} \in \mathcal{Y}^{strict}$. The
vertical 1-categories, horizontal 2-categories, and holonomies of
$\mathcal{M}\mathcal{L}\mathbb{D}$ and $\mathbb{D}$ are in fact the
same, and the squares correspond under the bijections
$\Lambda^{f,k}_{j,g}$.

Lastly, we verify that $\mathcal{M}$ is locally an isomorphism.
Clearly, $$\xymatrix@1{\mathcal{M}_{P,P'}:\mathcal{Z}^{strict}(P,P')
\ar[r] & \mathcal{Y}^{strict}(\mathcal{M}P,\mathcal{M}P')}$$ is
injective on objects and locally injective. If $F \in
\mathcal{Y}^{strict}(\mathcal{M}P,\mathcal{M}P')$, then
$\mathcal{M}\mathcal{L}F=F$, and similarly for the morphisms in
$\mathcal{Y}^{strict}(\mathcal{M}P,\mathcal{M}P')$.
\end{pf}

\begin{rmk} \label{brownmosaspecialcase}
Theorem \ref{YZ} is an $I$-category analogue of the equivalence in
\cite{brownmosa99} and \cite{spencer77} between the category of edge-symmetric double
categories with thin structure and the category of small
2-categories.
\end{rmk}

\begin{thm} \label{XY}
Let $I$ be a groupoid. The 2-category $\mathcal{X}^{strict}$ of
$I$-catego-ries (strict 2-algebras over the 2-theory of categories
with underlying groupoid $I$) is 2-equivalent to the 2-category
$\mathcal{Y}^{strict}$ of double categories $\mathbb{D}$ with
folding such that $(\mathbf{V}\mathbb{D})_0=I$ as defined in
Notation \ref{Xnotation} and Notation \ref{Ynotation}.

\end{thm}
\begin{pf}
We construct a 2-equivalence
$\xymatrix@1{\mathcal{J}:\mathcal{X}^{strict} \ar[r] &
\mathcal{Y}^{strict}}$. Suppose $X$ is an object of
$\mathcal{X}^{strict}$. From the strict 2-functor $\xymatrix@1{X:I^2
\ar[r] & Cat}$ we define $(\mathbf{V} \mathcal{J}(X))_0:=I$ and
$\mathbf{H}\mathcal{J}(X)(A,B):=X_{A,B}$. For $j$ as in
(\ref{sourcetarget}) the holonomy is
$\overline{j}:=X_{1^v_A,j}(1^h_A)$, which is the same as $P$ in
Lemma \ref{holonomyconstruction}. The squares
$$\xymatrix@R=3pc@C=3pc{A \ar[d]_j \ar[r]^{f} \ar@{}[dr]|{\alpha}
 & B \ar[d]^{k} \\ C \ar[r]_{g} & D }$$ are the
morphisms from $\overline{k} \circ f$ to $g \circ \overline{j}$ in
$X_{A,D}$, so that the bijection $\Lambda^{f,k}_{j,g}$ of Definition
\ref{foldingstructure} is the identity map. The horizontal and
vertical compositions of squares are defined by the axioms for
$\Lambda$ in Definition \ref{foldingstructure}.

It follows from the definitions that $\mathbf{H}\mathcal{J}(X)$ and
$\mathbf{V}\mathcal{J}(X)$ are 2-catego-ries. The associativity
axioms, identity axioms, and interchange law axiom for composition
of squares of $\mathcal{J}(X)$ follow from the analogous axioms for
the underlying 2-category $\mathbf{H}\mathcal{J}(X)$ of $X$ by
Definition \ref{foldingstructure}. In fact, the verifications are
formally similar to the analogous verifications for the double
category of quintets of a 2-category. Thus $\mathcal{J}(X)$ is a
double category with folding and belongs to $\mathcal{Y}^{strict}$.

The strict 2-functor $\mathcal{J}$ is defined similarly on
morphisms. For any morphism $\xymatrix@1{F:X \ar[r] & Y}$ in
$\mathcal{X}^{strict}$, then double functor $\mathcal{J}(F)$ is
defined as the identity on $I$. On $\mathbf{H}\mathcal{J}(X)(A,B)$
it is $F_{A,B}$, which extends to a definition on squares via
$\Lambda$. A naturality argument shows
$\mathcal{J}(F)(\overline{k})=\overline{\mathcal{J}(F)(k)}$, so that
$\mathcal{J}(F)$ is a morphism in $\mathcal{Y}^{strict}$.

If $\xymatrix@1{\sigma:F \ar@{=>}[r] & G}$ is a 2-cell in
$\mathcal{X}^{strict}$, then $\mathcal{J}(\sigma)_{A,B}$ is simply
$\sigma_{A,B}$. Since $\sigma$ is a modification, we have
$$Y_{j,k}(\sigma_{A,B}^f)=\sigma_{C,D}^{X_{j,k}(f)},$$ which implies
$$\sigma^{\overline{j}}_{A,C}=i_{\overline{j}}^v$$ by Lemma
\ref{holonomyconstruction}.  Furthermore, $\mathcal{J}(\sigma)$ is
compatible with horizontal composition and square identity because
$\sigma$ is. This concludes the definition of the strict 2-functor
$\mathcal{J}$.

We claim that $\mathcal{J}$ is surjective on objects up to
isomorphism. Let $\mathbb{D}$ be an object of
$\mathcal{Y}^{strict}$. Then $\mathbf{V}\mathbb{D}(A,B)=I(A,B)$. For
$A,B \in Obj \hspace{1mm} I$, let $X_{A,B}$ be the category whose
objects are horizontal morphisms $\xymatrix@1{f:A \ar[r] & B}$ in
$\mathbb{D}$ and whose morphisms are the squares of $\mathbb{D}$
with left and right vertical morphisms $1^v_A$ and $1^v_B$
respectively. For $f\in X_{A,B}$, a square $\alpha$ in $X_{A,B}$,
and vertical morphisms $\xymatrix@1{j:A \ar[r] & C}$ and
$\xymatrix@1{k:B \ar[r] & D}$, define
$$X_{j,k}(f)=[\overline{j^{-1}} \ \ f \ \ \overline{k}] \hspace{.25in}
X_{j,k}(\alpha):=[i_{\overline{j^{-1}}} \ \ \alpha \ \
i_{\overline{k}}].$$ Then $\xymatrix@1{X:I^2 \ar[r] & Cat}$ is a
clearly a 2-functor by the properties of $\mathbb{D}$, and even an
$I$-category. Moreover, $\mathbb{D}$ is isomorphic to
$\mathcal{J}(X)$ in $\mathcal{Y}^{strict}$ (squares of $\mathbb{D}$
are mapped to squares of $\mathcal{J}(X)$ using
$\Lambda^{\mathbb{D}}$).

The 2-functor $\mathcal{J}$ is locally an isomorphism by inspection.
Hence $\mathcal{J}$ is a 2-equivalence.
\end{pf}

Next we compare $\mathcal{X}^{strict}$ and $\mathcal{Z}^{strict}$.
The result is a strict version of Theorem 5.2 with trivial $T$ in
\cite{fiore3} improved from biequivalence to 2-equivalence in
Theorem \ref{XZ}. It is a corollary of Theorem \ref{YZ} and Theorem
\ref{XY}, but we present a direct proof:

\begin{thm} \label{XZ}
Let $I$ be a groupoid. The 2-category $\mathcal{X}^{strict}$ of
$I$-catego-ries (strict 2-algebras over the 2-theory of categories
with underlying groupoid $I$) is 2-equivalent to the 2-category
$\mathcal{Z}^{strict}$ of strict 2-functors $\xymatrix@1{P:I \ar[r]
& \mathcal{C}}$ that are the identity on objects as in Notation
\ref{Xnotation} and Notation \ref{Znotation}.
\end{thm}
\begin{pf}
We construct a 2-equivalence $\xymatrix@1{\mathcal{K}:\mathcal{X}^{strict} \ar[r] &
\mathcal{Z}^{strict}}$.
For an object
$\xymatrix@1{X:I^2 \ar[r] & Cat}$ of $\mathcal{X}$, we obtain a
strict 2-functor
$$\xymatrix{\mathcal{K}(X)=P:I \ar[r] & \mathcal{C}}$$
that is the identity on objects as in Lemma
\ref{holonomyconstruction}.

We define $\mathcal{K}$ compatibly on morphisms and 2-cells. Let
$\xymatrix@1{F:X \ar[r] & X'}$ be a morphism in
$\mathcal{X}^{strict}$, \ie $F$ is a strict 2-natural transformation
from $X$ to $X'$ which preserves composition and identity. We define
the 2-functor $\xymatrix@1{\mathcal{K}(F):\mathcal{C} \ar[r] &
\mathcal{C}'}$ to be the identity on objects, and as
$\xymatrix@1{F_{A,B}:X_{A,B} \ar[r] & X'_{A,B}}$ on
$Mor_{\mathcal{C}}(A,B)$ =$X_{A,B}$. Then $P' =\mathcal{K}(F) \circ
P$ because $F$ is 2-natural and preserves identity morphisms. For a
2-cell $\xymatrix@1{\sigma:F \ar@{=>}[r] & G}$ in
$\mathcal{X}^{strict}$, the natural transformation
$$\xymatrix@1{\mathcal{K}(\sigma)_{A,B}:\mathcal{K}(F)_{A,B} \ar[r] &
\mathcal{K}(G)_{A,B}}$$ is simply $\sigma_{A,B}$. One can easily
check that $\mathcal{K}$ is a strict 2-functor.

The strict 2-functor $\mathcal{K}$ is surjective on objects. If
$\xymatrix@1{P:I \ar[r] & \mathcal{C}}$ is an object of
$\mathcal{Z}^{strict}$, then $X_{A,B}:=Mor_{\mathcal{C}}(A,B)$ and
$X_{j,k}(f):=P(k) \circ f \circ P(j^{-1})$ defines an object of
$\mathcal{X}^{strict}$ which maps to $P$.

The strict 2-functor $\mathcal{K}$ is locally an isomorphism of
categories. It is clearly injective on morphisms and 2-cells. If
$\overline{F}$ is a morphism in $\mathcal{Z}^{strict}$ from
$\mathcal{K}(X)$ to $\mathcal{K}(X')$, then a pre-image is
necessarily defined by $F_{A,B}:=\overline{F}_{A,B}$, the
2-naturality of which is easily verified:
$$\aligned \overline{F}X_{j,k}(f)& =F(P(k) \circ f \circ P(j^{-1})) \\
&=FP(k) \circ F(f) \circ FP(j^{-1}) \\ &=P'(k) \circ F(f) \circ
P'(j^{-1}) \\ &= X'_{j,k}\overline{F}(f).
\endaligned$$
If $\xymatrix@1{\overline{\sigma}:\overline{F} \ar@{=>}[r] &
\overline{G}}$ is a 2-cell in $\mathcal{Z}^{strict}$, then a
pre-image is defined by $\sigma_{A,B}:=\overline{\sigma}_{A,B}$.
Since $\overline{\sigma}_{A,C}^{P(j)}=i_{P'(j)}$, we know that
$\sigma$ is a modification by Lemma \ref{holonomyconstruction}. The
modification $\sigma$ is clearly compatible with composition and
identity because $\overline{\sigma}$ is.

We conclude that $\mathcal{K}$ is a 2-equivalence of 2-categories.
\end{pf}

\begin{rmk}
A strict 2-functor $\xymatrix@1{P:I \ar[r] & \mathcal{C}}$ from a
1-category $I$ to a 2-category $\mathcal{C}$ (with the same object
set) that is the identity on objects is a special case of the notion
{\it weak equipment} in \cite{veritythesis}. There Verity constructs
a {\it double bicategory of squares} from a weak equipment, which
essentially defines our 2-functor
$\xymatrix@1{\mathcal{M}:\mathcal{Z}^{strict} \ar[r] &
\mathcal{Y}^{strict}}$ and a 2-functor $\xymatrix@1{\mathcal{Z}
\ar[r] & \mathcal{Y}}$ between the 2-categories defined in Section
\ref{pseudosection}, though the 2-cells are different in present
paper. Here we have constructed 2-equivalences between $\mathcal{X}$
and $\mathcal{Y}$ as well as between $\mathcal{X}$ and $\mathcal{Z}$
using connection pairs and foldings (in the strict and pseudo
cases).
\end{rmk}

In this section we proved the strict version of our desired result:
an $I$-category can be viewed as a double category with folding or
as a 2-functor from a 1-category into a 2-category. Since foldings
are equivalent to connection pairs, and edge-symmetric double groups
with connection pair\footnote{Whenever an edge-symmetric double
category is equipped with a connection pair, we assume the holonomy
to be trivial.} are equivalent to crossed modules, one can expect
that Theorem \ref{XY} has implications for crossed modules. Indeed,
we pursue this in the next section.

\section{2-Groups, Double Groups, and Crossed Modules}
It is often useful to investigate one-object cases of categorical
concepts and compare them with more familiar concepts. For example,
a one-object category is simply a monoid, and a one-object groupoid
is a group. In this section, we investigate one-object
$I$-categories with everything invertible and compare them with
other notions in the literature: crossed modules and double groups.

To see how our comparison of one-object-everything-invertible
$I$-categories will work, consider 2-groupoids. These are
2-categories in which every morphism and every 2-cell is invertible.
We call a one-object 2-groupoid a {\it 2-group}.\footnote{In this
article all 2-groupoids are strict.} This is the same as a group
object in $Cat$, or {\it categorical group} as in Theorem
\ref{categoricalgroups=2groups}. Though the notion of 2-group is no
more familiar than the notion of 2-groupoid, we can compare it to
something familiar. Brown and Spencer proved in
\cite{brownspencer74} (and attribute the result to Verdier) that
categorical groups (and hence 2-groups) are equivalent to {\it
crossed modules}. This last concept is much more familiar to
topologists than 2-groups. Whitehead first introduced crossed
modules in \cite{whitehead49} and proved with Mac Lane that they
model path-connected homotopy 2-types in \cite{maclanewhitehead50}.
The survey \cite{paolisurvey} contains an account of the use of
crossed modules and their higher-dimensional analogues to model
homotopy types. Recently, 2-groups have been studied in
\cite{baezlauda}. Our comparison of one-object $I$-categories with
everything invertible will build on this result of Brown and
Spencer. In fact, Brown and Spencer obtained a 2-equivalence, and we
will in Theorem \ref{WZ} as well.

In addition to the comparison with crossed modules, we also compare
with double groups. A {\it double groupoid} is a double category in
which all morphisms and squares are iso. In particular, squares are
required to be isos under both vertical composition and horizontal
composition. In analogy to 2-groups, we shall call a one-object
double groupoid a {\it double group}\footnote{A group object in the
category of groups is simply an abelian group by Eckmann-Hilton.
Thus double groups and group objects in the category of groups are
not the same.}. Brown and Spencer proved that edge-symmetric double
groups with connection\footnote{The holonomies here are trivial.}
are equivalent to crossed modules in \cite{brownspencer76}. We
extend this to a 2-equivalence between general double groups with
folding and crossed modules under groups in Theorem
\ref{specialextension}. Brown and Higgins showed in
\cite{brownhigginscubes} that so-called crossed modules over
groupoids ({\it not} in the sense of an over category) are
equivalent to edge-symmetric double groupoids with connection. In
\cite{brownmackenzie}, Brown and Mackenzie substantially generalized
\cite{brownspencer76} to an equivalence between locally trivial (not
necessarily edge-symmetric) double Lie groupoids and so-called
locally trivial core diagrams. This is an equivalence between double
groupoids with certain filling conditions and core diagrams. Their
Theorem 4.2 treats the discrete case as well. Double groupoids have
recently found application in the theory of weak Hopf algebras in
\cite{andruskiewitschnatalequantum} and
\cite{andruskiewitschnataletensor}.

In most cases, our 2-cells are the vertical natural transformations.
Double categories provide a good context for the 2-equivalence of
categorical groups, 2-groups, and crossed modules. The 2-natural
transformations of 2-groups do {\it not} correspond to the {\it
homotopies} in the 2-category of crossed modules, instead one needs
the vertical transformations. Theorem \ref{2groupequivcrossedmodule}
and Theorem \ref{specialextension} also hold for the horizontal
natural transformations after adjusting the notion of 2-cell in the
2-category of crossed modules.

We begin by stating Theorems \ref{YZ}, \ref{XY}, and \ref{XZ} in the
special case that $I$ is a group and all morphisms in each structure
are invertible. We call these sub-2-categories $\mathcal{X}^{inv},
\mathcal{Y}^{inv},$ and $\mathcal{Z}^{inv}$ in Notation
\ref{oneobjectNotation}. After recalling 2-groups, categorical
groups, crossed modules, and the proof by Brown and Spencer, we show
that $\mathcal{Z}^{inv}$ is 2-equivalent to the 2-category
$\mathcal{W}^{inv}$ of crossed modules under $\xymatrix@1{\{e\}
\ar[r] & I}$. Therefore $I$-categories which have only one object
and everything invertible are 2-equivalent to crossed modules under
$I$. The 2-equivalence of $\mathcal{Y}^{inv}$ and
$\mathcal{W}^{inv}$ says that the 2-category of double groups
$\mathbb{D}$ with folding such that $(\mathbf{V}\mathbb{D})_0=I$ is
2-equivalent to the 2-category of crossed modules under $I$. Lastly
we turn to double groups with folding.

\begin{notation} \label{oneobjectNotation}
In this section, $I$ denotes a one-object groupoid, \ie a group. Let
$\mathcal{X}^{inv}$ denote the 2-category of $I$-categories $X$ with
$X_{*,*}$ a groupoid whose objects and morphisms are invertible with
respect to the 2-algebra composition $\circ$. Let
$\mathcal{Y}^{inv}$ denote the 2-category of double groups
$\mathbb{D}$ with folding such that $(\mathbf{V}\mathbb{D})_0=I$.
Let $\mathcal{Z}^{inv}$ denote the 2-category of strict 2-groups
$\mathcal{C}$ equipped with a strict 2-functor $\xymatrix@1{P:I
\ar[r] & \mathcal{C}}$. Morphisms and 2-cells of the 2-categories
$\mathcal{X}^{inv},\mathcal{Y}^{inv},$ and $\mathcal{Z}^{inv}$ are
as in the respective categories of Notations \ref{Xnotation},
\ref{Ynotation}, and \ref{Znotation}.
\end{notation}

\begin{thm} \label{XYZinv}
The 2-categories $\mathcal{X}^{inv},\mathcal{Y}^{inv},$ and
$\mathcal{Z}^{inv}$ are 2-equivalent.
\end{thm}
\begin{pf}
The 2-equivalence of $\mathcal{X}^{inv},\mathcal{Y}^{inv},$ and
$\mathcal{Z}^{inv}$ follows from Theorems \ref{YZ}, \ref{XY}, and
\ref{XZ}.
\end{pf}

\begin{defn}
A one-object 2-category in which all 1-cells and all 2-cells are
invertible is called a {\it 2-group}. We view a 2-group
$\mathcal{C}$ as a double category with one object, trivial vertical
morphisms,  and with horizontal morphisms and squares given by the
1-cells and 2-cells of $\mathcal{C}$. In other words, we view a
2-group $\mathcal{C}$ as $\mathbb{H}\mathcal{C}$. A morphism of
2-groups is a 2-functor. This is the same as a double functor
between the associated double categories. A {\it 2-cell} is a
vertical natural transformation, {\it not} a 2-natural
transformation. We denote this 2-category by $2$-$Gp$. It is a
sub-2-category of $Cat(Cat)_v$, the 2-category of double categories,
double functors, and vertical natural transformations.
\end{defn}

\begin{defn}
A {\it categorical group} is a group object in $Cat$. This is a
category $(X_0, X_1)$ equipped with a functor
$$\xymatrix{(X_0,X_1) \times (X_0,X_1) \ar[r] & (X_0,X_1)}$$ which strictly satisfies the
axioms of a group. A {\it morphism of categorical groups} is a
functor compatible with group structure. A {\it 2-cell} is a natural
transformation compatible with group structure. We denote the
2-category of group objects in $Cat$ by $Gp(Cat)$.
\end{defn}

\begin{thm} \label{categoricalgroups=2groups}
The 2-category of categorical groups, morphisms, and 2-cells is
2-equivalent to the 2-category of 2-groups, 2-functors, and vertical
natural transformations.
\end{thm}
\begin{pf}
The inclusion of $Gp$ into $Cat$ induces an inclusion of
2-categories $\xymatrix@1{Gp(Cat) \ar[r] & Cat(Cat)_v}$. This
assigns a categorical group $(X_0,X_1)$ to the double category with
one object, no nontrivial vertical morphisms, horizontal morphisms
$X_0$, and squares $X_1$. The horizontal composition is the group
operation, and the vertical composition of squares is composition in
the category $(X_0,X_1)$. Morphisms and 2-cells of $Gp(Cat)$ are
mapped to double functors and {\it vertical} natural
transformations. Thus $Gp(Cat)$ is contained in $2$-$Gp$. Every
2-group is isomorphic to one in $Gp(Cat)$, so $Gp(Cat)$ is
2-equivalent to $2$-$Gp$.
\end{pf}

\begin{defn}
A {\it crossed module} $\xymatrix@1{\partial:H \ar[r] & G}$ consists
of
\begin{itemize}
\item
groups $H$ and $G$
\item
a group homomorphism $\xymatrix@1{\partial:H \ar[r] & G}$
\item
a left action of $G$ on $H$ by automorphisms written $(g,\alpha)
\mapsto {}^g\alpha$ such that:
\begin{enumerate}
\item
$\partial({}^g\alpha)=g\partial(\alpha)g^{-1}$ for all $\alpha \in
H$ and $g \in G$,
\item
${}^{\partial (\beta)} \alpha=\beta \alpha \beta^{-1}$ for all
$\alpha,\beta \in H$.
\end{enumerate}
\end{itemize}
\end{defn}
\begin{defn}
If $(H,G,\partial)$ and $(H',G',\partial')$ are crossed modules,
then a {\it morphism} $\xymatrix@1{(p,q):(H,G,\partial) \ar[r] &
(H',G',\partial')}$ consists of group homomorphisms $p$ and $q$ such
that the following diagram commutes
$$\xymatrix{H \ar[r]^{\partial} \ar[d]_p & G \ar[d]^q \\ H' \ar[r]_{\partial'} & G'}$$
and $p({}^g\alpha)={}^{q(g)}p(\alpha)$ for all $g \in G$ and $\alpha
\in H$.
\end{defn}

\begin{defn} \label{crossedmodulehomotopy}
If $\xymatrix@1{(p_1,q_1), (p_2,q_2):(H,G,\partial) \ar[r] &
(H',G',\partial')}$ are morphisms of crossed modules, then a {\it
homotopy} $\xymatrix@1{\nu:(p_1,q_1) \ar@{=>}[r] & (p_2,q_2)}$ is a
function $\xymatrix@1{\nu:G \ar[r] & H'}$ such that
$(\partial'\nu(f))q_1(f)=q_2(f)$ and:
\begin{enumerate}
\item
For all $f, g \in G$ and $\alpha \in H$ such that
$\partial(\alpha)f=g$, we have
$$p_2(\alpha)\nu(f)=\nu(g)p_1(\alpha),$$
\item
For all $f,g \in G$, the {\it derivation rule} holds:
$$\nu(g)\cdot {}^{q_1(g)}\nu(f)=\nu(gf).$$
\end{enumerate}
The vertical composition of homotopies
$$\xymatrix{(p_1,q_1) \ar@{=>}[r]^{\nu_1} & (p_2,q_2) \ar@{=>}[r]^{\nu_2} & (p_3,q_3)}$$ is $\xymatrix{f \ar@{|->}[r] &
\nu_2(f)\nu_1(f).}$ The horizontal composition of homotopies
$$\xymatrix{(H_1,G_1,\partial_1) \ar@/^2pc/[rr]_{\quad}^{(p_1,q_1)}="1" \ar@/_2pc/[rr]_{(p_1',q_1')}
="2" & & \ar@{}"1";"2"|(.2){\,}="7" \ar@{}"1";"2"|(.8){\,}="8"
\ar@{=>}"7" ;"8"^{\nu_1} (H_2,G_2,\partial_2)
\ar@/^2pc/[rr]_{\quad}^{(p_2,q_2)}="1"
\ar@/_2pc/[rr]_{(p_2',q_2')}="2" & & \ar@{}"1";"2"|(.2){\,}="7"
\ar@{}"1";"2"|(.8){\,}="8" \ar@{=>}"7" ;"8"^{\nu_2}
(H_3,G_3,\partial_3)}$$ is $\xymatrix{f \ar@{|->}[r] &
\nu_2(q_1'(f))\cdot p_2(\nu_1(f)).}$
\end{defn}

Crossed modules, morphisms, and homotopies form a 2-category denoted
$XMod$. For more on crossed modules as internal categories and their
2-cells, see \cite{datuashvili}. Homotopies and derivations for more general
crossed modules as needed for a 2-dimensional notion of
holonomy are considered in \cite{brownicenhomotopies}.

\begin{examp}
 An example of a crossed module is the inclusion of a normal
subgroup $H$ into a group $G$ where the action is conjugation by
elements of $G$. In particular, $\xymatrix@1{\{e\} \ar[r] & I}$ is a
crossed module for any group $I$. We abbreviate $\xymatrix@1{\{e\}
\ar[r] & I}$ by $I$.
\end{examp}

\begin{examp}[Whitehead]
Let $(X,A,*)$ be a pair of based spaces. Then the boundary map
$\xymatrix@1{\partial:\pi_2(X,A,*) \ar[r] & \pi_1(A,*)}$ is a
crossed module with action given by the standard action of the
fundamental group. Crossed modules are known to model pointed path-connected weak homotopy
2-types algebraically. A proof is sketched in \cite{brown99}.
\end{examp}

In preparation for our theorem, we summarize Brown and Spencer's
proof as recounted in \cite{forrester}. Brown and Spencer originally
showed that categorical groups are 2-equivalent to crossed modules,
crossed module morphisms, and homotopies. The 2-category of
categorical groups is 2-equivalent to the 2-category of 2-groups,
functors, and vertical natural transformations. Since we are
interested in double groups, we work with the latter 2-category of
2-groups instead of categorical groups. See \cite{brownhigginstensor} for the analogue of
Theorem \ref{2groupequivcrossedmodule} in arbitrary dimensions.

\begin{thm}[Brown-Spencer in \cite{brownspencer74}] \label{2groupequivcrossedmodule}
The 2-category $2$-$Gp$ of 2-groups, functors, and vertical natural
transformations is 2-equivalent to the 2-category $XMod$ of crossed
modules, crossed module morphisms, and homotopies.
\end{thm}
\begin{pf}
Let $\mathcal{C}$ be a 2-group. We obtain a crossed module from
$\mathcal{C}$ as follows. The group $G$ consists of the objects of
$Mor_{\mathcal{C}}(*,*)$. In particular, $e_G$ is the identity
morphism. The group $H$ consists of 2-cells $\alpha$ in
$\mathcal{C}$ whose source is $e_G$ and $\partial$ is the target
map.
$$
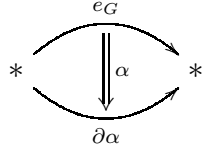
\xymatrix{ \ast \ar@/^1.5pc/[rr]_{\quad}^{e_G}="1"
\ar@/_1.5pc/[rr]_{\partial \alpha}="2" & & \ast
\ar@{}"1";"2"|(.2){\,}="7" \ar@{}"1";"2"|(.8){\,}="8" \ar@{=>}"7"
;"8"^{\alpha} }$$ The group $G$ acts on $H$ on the left by
conjugation, in other words ${}^g\alpha$ has the form below.
$$
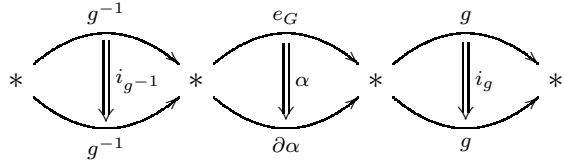
\xymatrix{\ast \ar@/^1.5pc/[rr]_{\quad}^{g^{-1}}="1" \ar@/_1.5pc/[rr]_{g^{-1}}
="2" & & \ar@{}"1";"2"|(.2){\,}="7" \ar@{}"1";"2"|(.8){\,}="8"
\ar@{=>}"7" ;"8"^{i_{g^{-1}}} \ast
\ar@/^1.5pc/[rr]_{\quad}^{e_G}="1" \ar@/_1.5pc/[rr]_{\partial
\alpha}="2" & & \ar@{}"1";"2"|(.2){\,}="7"
\ar@{}"1";"2"|(.8){\,}="8" \ar@{=>}"7" ;"8"^{\alpha} \ast
\ar@/^1.5pc/[rr]_{\quad}^{g}="1" \ar@/_1.5pc/[rr]_{g} ="2" & &
\ar@{}"1";"2"|(.2){\,}="7" \ar@{}"1";"2"|(.8){\,}="8" \ar@{=>}"7"
;"8"^{i_{g}} \ast}$$

If $\xymatrix@1{F:\mathcal{C} \ar[r] & \mathcal{C}'}$ is a
2-functor, then we obtain a morphism of crossed modules by
restricting $F$ to $G$ and $H$. A 2-cell $\xymatrix@1{F_1
\ar@{=>}[r] & F_2}$ in the 2-category of 2-groups is a natural
transformation
$$\xymatrix@1{\sigma:F_1|_{Mor_{\mathcal{C}}(*,*)} \ar@{=>}[r] &
F_2|_{Mor_{\mathcal{C}'}(*,*)}}$$ that is compatible with the
horizontal composition of 1- and 2-cells. We obtain a homotopy
$\xymatrix@1{\nu:G \ar[r] & H'}$ by defining
$\nu(g):=\sigma^gi_{(F_1g)^{-1}}$. Here concatenation denotes the
horizontal composition of 2-cells. The naturality corresponds to (i)
and the compatibility with horizontal composition corresponds to
(ii) in Definition \ref{crossedmodulehomotopy}.

Next we describe how to get a 2-group $\mathcal{C}$ from a crossed
module $\xymatrix@1{\partial:H \ar[r] & G}$. The set of morphisms of
$\mathcal{C}$ is $G$ and the set of 2-cells of $\mathcal{C}$ is $H
\rtimes G$. The source and target of the 2-cell $(\alpha, g)$ are
$g$ and $\partial (\alpha) g$ respectively. Horizontal composition
of 2-cells is given by the group operation in $H \rtimes G$ and
vertical composition is
$$(\alpha_2,\partial(\alpha_1)g_1) \odot (\alpha_1, g_1):=(\alpha_2\alpha_1,
g_1).$$ The vertical identity is $\xymatrix@1{i_g:=(e_H,g):g
\ar@{=>}[r] & g}$.

If $(p,q)$ is a morphism of crossed modules, we obtain a 2-functor
$\xymatrix@1{\mathcal{C} \ar[r] & \mathcal{C}'}$ as $q$ on 1-cells
and $(p,q)$ on 2-cells. A homotopy $\nu$ in the 2-category of
crossed modules gives rise to a natural transformation
$\xymatrix@1{\sigma:F_1|_{Mor_{\mathcal{C}}(*,*)} \ar@{=>}[r] &
F_2|_{Mor_{\mathcal{C}'}(*,*)}}$ by defining
$\sigma^g:=(\nu(g),q_1(g))$. Further, this natural transformation is a vertical
natural transformation: the derivation rule for homotopies guarantees that
$\sigma$ is compatible with composition of horizontal morphisms as in Definition
\ref{verticalnaturaltransformation} (ii), since
$$\aligned \sigma^g \sigma^f &=(\nu(g),q_1(g))(\nu(f),q_1(f)) \\
&= (\nu(g)\cdot {}^{q_1(g)}\nu(f),q_1(g)q_1(f)) \\
&= (\nu(gf),q_1(gf)) \\
&= \sigma^{gf}.
\endaligned$$

The composite 2-functor from crossed modules to 2-groups and to
crossed modules back again is 2-naturally isomorphic to the
identity. On the other hand, if we start with a 2-group
$\mathcal{C}$, and take the associated crossed module, note that the
group of 2-cells of $\mathcal{C}$ (under horizontal composition) is
isomorphic to $H \rtimes G$ by the map which sends
$$\xymatrix{ \ast \ar@/^1.5pc/[rr]_{\quad}^{g_1}="1"
\ar@/_1.5pc/[rr]_{g_2}="2" & & \ast \ar@{}"1";"2"|(.2){\,}="7"
\ar@{}"1";"2"|(.8){\,}="8" \ar@{=>}"7" ;"8"^{\gamma} }$$ to $(\gamma
i_{g_1^{-1}},g_1)$. Using this map, one can see that the composite
2-functor from 2-groups to crossed modules and to 2-groups back
again is 2-naturally isomorphic to the identity.
\end{pf}

\begin{notation}
The objects of the 2-category $\mathcal{W}^{inv}$ are crossed
modules under $I$. These are crossed modules $\xymatrix@1{\partial:H
\ar[r] & G}$ equipped with a morphism of crossed modules
$$\xymatrix{\{e\} \ar[r] \ar[d] & I \ar[d]^P \\ H \ar[r]_{\partial} & G.}$$
Morphisms of $\mathcal{W}^{inv}$ are morphisms $(p,q)$ of crossed
modules under $I$, in other words
$$\xymatrix{I \ar[r]^P \ar[dr]_{P'} & G \ar[d]^q \\ & G'}$$ commutes.
A 2-cell in $\mathcal{W}^{inv}$ is a homotopy $\nu$ such that
$\nu(P(j))=e_{H'}$ for all $j \in I$.
\end{notation}

\begin{thm} \label{WZ}
The 2-category $\mathcal{W}^{inv}$ of crossed modules under $I$ is
2-equivalent to the 2-category $\mathcal{Z}^{inv}$ of 2-groups under
$I$.
\end{thm}
\begin{pf}
The 2-equivalence from crossed modules to 2-groups in Theorem
\ref{2groupequivcrossedmodule} extends to a 2-equivalence
$\xymatrix@1{\mathcal{N}:\mathcal{W}^{inv} \ar[r] &
\mathcal{Z}^{inv}}$. A strict 2-functor $\xymatrix@1{I \ar[r] &
\mathcal{C}}$ is the same as a morphism of crossed modules from $I$
into the crossed module associated to $\mathcal{C}$. Morphisms of
crossed modules under $I$ are the same as morphisms of 2-groups
under $I$.

We observe that the 2-cells $\xymatrix@1{\nu:(p_1,q_1) \ar@{=>}[r] &
(p_2,q_2)}$ in $\mathcal{W}^{inv}$ are precisely the 2-cells
$\xymatrix@1{\mathcal{N}(p_1,q_1) \ar@{=>}[r] &
\mathcal{N}(p_2,q_2)}$ in $\mathcal{Z}^{inv}$. From Theorem
\ref{2groupequivcrossedmodule} we know that the homotopies
$\xymatrix@1{(p_1,q_1) \ar@{=>}[r] & (p_2,q_2)}$ in $XMod$
correspond to the 2-cells $\xymatrix@1{\mathcal{N}(p_1,q_1)
\ar@{=>}[r] & \mathcal{N}(p_2,q_2)}$ in $2$-$Gp$. It suffices to
show that $\nu(P(j))=e_{H'}$ if and only if its associated 2-group
2-cell $g \mapsto \sigma^g=(\nu(g), q_1(g))$ satisfies
$\sigma^{P(j)}=i_{P'(j)}$. But this is the case, since
$$\sigma^{P(j)}=(\nu(P(j)), q_1P(j))=(\nu(P(j)), P'(j))$$
$$i_{P'(j)}=(e_{H'},P'(j)).$$
Therefore $\mathcal{N}$ is a 2-functor that is essentially
surjective on objects and locally an isomorphism, \ie a
2-equivalence.

Alternatively, one could use the 2-equivalence
$\xymatrix@1{2\text{-}Gp \ar[r] & XMod}$ and similarly show that a
2-cell $\sigma$ in $2$-$Gp$ satisfies $\sigma^{P(j)}=i_{P'(j)}$ if
and only if the associated homotopy $\nu(g)=\sigma^gi_{(F_1g)^{-1}}$
satisfies $\nu(P(j))=e_{H'}$.
\end{pf}

With these notions we can extend Brown and Spencer's equivalence
between crossed modules and edge-symmetric double groups with
connection to the non-edge-symmetric case. The nontrivial holonomy
corresponds to a morphism of crossed modules from the vertical group
into the crossed module associated to the horizontal 2-group. In the
rest of this section, we no longer consider fixed $I$. Our proof
builds on the proof of Theorem \ref{2groupequivcrossedmodule}.

\begin{notation}
Let $Gp/XMod$ denote the 2-category of crossed modules under groups.
An object consists of a crossed module $(H,G,\partial)$ and a group
$I$ equipped with a crossed module morphism
$$\xymatrix{\{e\} \ar[r] \ar[d] & I \ar[d] \\ H \ar[r]_{\partial} & G.}$$
A morphism in $Gp/XMod$ is a morphism in the arrow category of
crossed modules. A 2-cell $\xymatrix@1{(r_1,p_1,q_1) \ar@{=>}[r] &
(r_2,p_2,q_2)}$ in $Gp/XMod$ is a homotopy
$\xymatrix@1{\nu:(p_1,q_1) \ar@{=>}[r] & (p_2,q_2)}$. Note that all
homotopies of crossed module morphisms $\xymatrix@1{I \ar[r] & I'}$
are trivial, so we do not include this in the data for a 2-cell in
$Gp/XMod$.

Let $DblGpFold$ denote the 2-category of double groups with folding.
The morphisms are morphisms of double categories with folding as in
Definition \ref{morphismwithfolding}. The 2-cells $\xymatrix@1{F_1
\ar@{=>}[r] & F_2}$ are {\it vertical} natural transformations
between the restrictions of $F_1$ and $F_2$ to the sub-double
category with trivial vertical morphisms. We do {\it not} require
that the vertical natural transformations are compatible with
folding.
\end{notation}

The 2-category $DblGpFold$ is like $\mathcal{Y}^{inv}$, except that
we allow $I$ to vary and do not require the 2-cells to be compatible
with folding. To extend the equivalence in \cite{brownspencer76}
between edge-symmetric double groups with connection and crossed
modules to a 2-equivalence, one is forced to take vertical
transformations with identity components as the 2-cells between
morphisms of edge-symmetric double groups (2-equivalences are local
isomorphisms). Likewise, in the non-edge-symmetric case of
$DblGpFold$, the vertical transformations are not required to be
compatible with folding: any vertical transformation in $DblGpFold$
with identity components that is compatible with folding (as in
Definition \ref{2cellwithfolding}) is necessarily trivial.

Our choice of 2-cell in $DblGpFold$ is compatible with the degree 1
part of the internal hom for cubical $\omega$-groupoids constructed
in \cite{brownhigginstensor}: an edge-symmetric double groupoid with
connection is a 2-truncated cubical $\omega$-groupoid as defined in
\cite{brownhigginscubes}. We now extend the equivalence in
\cite{brownspencer76} between edge-symmetric double groups with
connection and crossed modules to the non-edge-symmetric setting and
upgrade it to a 2-equivalence:

\begin{thm} \label{specialextension}
The 2-category $Gp/XMod$ of crossed modules under groups is
2-equivalent to the 2-category $DblGpFold$ of double groups with
folding, morphisms, and vertical transformations.
\end{thm}
\begin{pf}
We define a 2-equivalence $\xymatrix@1{\mathcal{R}:DblGpFold \ar[r]
& Gp/XMod}$. For a double group $\mathbb{C}$ with folding, we define
$I$ to be the group of vertical morphisms, $(H,G,\partial)$ to be
the crossed module associated to the horizontal 2-group, and the
homomorphism $\xymatrix@1{I \ar[r] & G}$ to be the holonomy, so that
$\mathcal{R}(\mathbb{C})=(I,H,G,\partial)$. If
$\xymatrix@1{F:\mathbb{C} \ar[r] & \mathbb{C}'}$ is a morphism, then
$\mathcal{R}(F)$ is the restriction of $F$ to $I,H$, and $G$. If
$\sigma$ is a vertical transformation $\xymatrix@1{F_1 \ar@{=>}[r] &
F_2}$ with identity components, then $\mathcal{R}(\sigma)$ is the
homotopy associated to the restriction of $\sigma$ to the horizontal
2-group as in Theorem \ref{2groupequivcrossedmodule}.

The 2-functor $\mathcal{R}$ is surjective on objects up to
isomorphism. If $(I,H,G,\partial)\in Gp/XMod$, then we construct the
2-group $\mathcal{C}$ associated to the crossed module
$(H,G,\partial)$ as in Theorem \ref{2groupequivcrossedmodule}. It
has horizontal morphisms $G$ and 2-cells $H \rtimes G$. The group
homomorphism $\xymatrix@1{I \ar[r] & G}$ determines a 2-functor
$\xymatrix@1{I \ar[r] & C}$. This data determines a double category
$\mathbb{C}$ with folding as in Theorem \ref{YZ}: the vertical
morphisms are $I$, the horizontal 2-category is the 2-group
$\mathcal{C}$, and the squares are determined by the 2-cells of the
horizontal 2-group by the folding. We see that
$\mathcal{R}(\mathbb{C}) \cong (I,H,G,\partial)$.

Lastly we verify that the functor
$$\xymatrix{\mathcal{R}_{\mathbb{C},\mathbb{C}'}:DblGpFold(\mathbb{C},\mathbb{C}')
\ar[r] & Gp/XMod(\mathcal{R}(\mathbb{C}),\mathcal{R}(\mathbb{C}')
)}$$ is an isomorphism of categories. Two morphisms
$\xymatrix@1{F_1,F_2:\mathbb{C} \ar[r] & \mathbb{C}'}$ that coincide
on the vertical 1-category and the horizontal 2-category also
coincide on the squares by the compatibility with folding.
Similarly, a morphism $\xymatrix@1{(r,p,q):\mathcal{R}(\mathbb{C})
\ar[r] & \mathcal{R}(\mathbb{C}')}$ has a pre-image because the
folding defines a morphism on general squares from the horizontal
2-functor associated to $(p,q)$ (as in Theorem
\ref{2groupequivcrossedmodule}) and the vertical 1-functor $r$. Thus
$\mathcal{R}_{\mathbb{C},\mathbb{C}'}$ is bijective on objects. The
2-cells $\xymatrix@1{\mathcal{R}(F_1) \ar@{=>}[r] &
\mathcal{R}(F_2)}$ are the 2-cells between the underlying
crossed-module morphisms of $\mathcal{R}(F_1)$ and
$\mathcal{R}(F_2)$. By Theorem \ref{2groupequivcrossedmodule}, the
latter are in bijective correspondence with the vertical
transformations between the restrictions of the double functors
$F_1$ and $F_2$ to the sub-double categories of $\mathbb{C}$ and
$\mathbb{C}'$ with trivial vertical morphisms. These are precisely
the 2-cells of $DblGpFold$. Hence
$\mathcal{R}_{\mathbb{C},\mathbb{C}'}$ is fully faithful and an
isomorphism of categories.

The 2-functor $\mathcal{R}$ is locally an isomorphism and surjective
up to isomorphism on objects, so that
$\xymatrix@1{\mathcal{R}:DblGpFold \ar[r] & Gp/XMod}$ is a
2-equivalence.

\end{pf}

\begin{rmk}
Theorem \ref{2groupequivcrossedmodule} used vertical transformations
as the 2-cells in  $2$-$Gp$ and homotopies as the 2-cells in $XMod$.
One could just as well work with horizontal transformations in
$2$-$Gp$ to obtain a 2-equivalence. However, the notion of 2-cell in
$XMod$ must be changed appropriately.  A 2-cell
$\xymatrix@1{w:(p_1,q_1) \ar@{=>}[r] & (p_2,q_2)}$ in this approach
is an element $w \in G'$ such that:
\begin{enumerate}
\item
$wq_1 (g)w^{-1}=q_2(g)$ for all $g \in G$,
\item
${}^wp(h)=p'(h)$ for all $h \in H$.
\end{enumerate}
If we use these 2-cells in $Gp/XMod$ and horizontal transformations
{\it compatible with folding} as 2-cells in $DblGpFold$, then we
obtain an analogue of Theorem \ref{specialextension}. The proof is
very similar, but relies on Remark \ref{2cellcompatiblewithfolding}
in the discussion of 2-cells. I thank Simona Paoli for pointing out
to me that the two notions of 2-cells correspond to horizontal and
vertical natural transformations.
\end{rmk}


This concludes our discussion of strict structures.

\section{Pseudo Double Categories with Folding} \label{pseudodoublesection}

Next we turn our attention to weak structures and work towards
pseudo versions of Theorems \ref{YZ}, \ref{XY}, and \ref{XZ}. There
are various ways of weakening a double category. Recall that a
double category contains two 2-categories, namely its horizontal and
vertical 2-categories as in Definition \ref{horizontal2category}.
One can weaken either or both of these to a bicategory, but in many
applications only one direction is typically weak. In this paper, we
prefer to make the horizontal 2-category into a horizontal
bicategory. This corresponds to the passage from category object in
$Cat$ to pseudo category object in $Cat$. Pseudo double categories,
and more generally pseudo category objects, have been studied in
\cite{grandisdouble1}, \cite{grandisdouble2},
\cite{martinspseudocategories}, and \cite{veritythesis}. Often one
can arrange the units of the horizontal bicategory to be strict.
Later we work with strict units.

\begin{defn} \label{pseudodoubledefn}
A {\it pseudo double category $\mathbb{D}$} consists of a class of
{\it objects}, a set of {\it horizontal morphisms}, a set of {\it
vertical morphisms}, and a set of {\it squares} with source and
target as in (\ref{sourcetarget}). The vertical morphisms are
equipped with a composition, as are the horizontal morphisms. The
squares are equipped with a vertical and a horizontal composition.
morphisms and squares also form a category under vertical
composition with identity squares $i^v_f$ as in
(\ref{verticalidentity}) which satisfy
$$[i_{f_1}^v \ \ i_{f_2}^v ] = i_{[f_1 \ f_2]}^v.$$ There are also
distinguished squares $i^h_j$ (not necessarily identity) as in (\ref{horizontalidentity}) which
satisfy
$$\begin{bmatrix} i_{j_1}^h \vspace{2mm} \\ i_{j_2}^h \end{bmatrix}=i_{\vcomp{j_1}{j_2 \vspace{1mm}}}^h
\quad \text{and} \quad i^h_{1^v_A}=i^v_{1_A^h}.$$ The objects,
horizontal morphisms, and squares with trivial left and right sides
form a bicategory with coherence iso 2-cells
$$\xymatrix@C=3pc@R=3pc{A \ar[r]^{1_B \circ f} \ar[d]_{1^v_A} \ar@{}[dr]|{\lambda_f}
& B \ar[d]^{1^v_B} & A \ar[r]^{f\circ 1_A} \ar[d]_{1^v_A}
\ar@{}[dr]|{\rho_f} & B \ar[d]^{1^v_B} & A \ar[r]^{h\circ (g \circ
f)} \ar[d]_{1^v_A} \ar@{}[dr]|{\alpha_{h,g,f}} & C \ar[d]^{1^v_C}
\\ A \ar[r]_{f} & B & A \ar[r]_{f} & B & A \ar[r]_{(h \circ g) \circ f} & C }$$
that satisfy the usual coherence triangle diagram and coherence
pentagon diagram for bicategories as in the original \cite{benabou1}, or in the
review \cite{leinsterbicat}, or in the Appendix to \cite{fiore3}.
These coherence iso 2-cells are also natural for all squares, for example
$$\begin{array}{c}
\xymatrix@R=3pc@C=3pc{A \ar[r]^{1^h_A} \ar@{}[dr]|{i^h_j} \ar[d]_j &
A \ar[d]|{\tb{j}} \ar[r]^f \ar@{}[dr]|\beta & B \ar[d]^k & A
\ar[r]^{1^h_A} \ar[d]_{1^v_A} \ar@{}[drr]|{\rho_f} & A \ar[r]^f & B
\ar[d]^{1^v_B}
\\ C
\ar[r]_{1^h_C} \ar[d]_{1^v_C} \ar@{}[drr]|{\rho_g} & C \ar[r]_g & D
\ar[d]^{1^v_D} \ar@{}[r]|= & A \ar[d]_j \ar@{}[rrd]|{\beta}
\ar[rr]|{\lr{f}} & & B \ar[d]^k \\ C \ar[rr]_g & & D &  C \ar[rr]_g
& & D}
\end{array}$$
and also for $\begin{array}{c} \xymatrix@R=3pc@C=3pc{ \ar[r]^{f_1}
\ar[d] \ar@{}[dr]|{\beta} &  \ar[r]^{f_2} \ar[d] \ar@{}[dr]|{\gamma}
& \ar[d] \ar[r]^{f_3} \ar@{}[dr]|{\delta} &  \ar[d] \\  \ar[r]_{g_1}
& \ar[r]_{g_2} & \ar[r]_{g_3} & }
\end{array}$ we have
$$\begin{array}{c}
\begin{bmatrix}
\ [ \ \beta \ \gamma \ ] \ \delta \ \\ \alpha_{g_3,g_2,g_1}
\end{bmatrix}
\end{array}=\begin{array}{c}
\begin{bmatrix}
\alpha_{f_3,f_2,f_1} \\
\ \beta \ [ \ \gamma \ \delta \ ] \
\end{bmatrix}.
\end{array}$$

Lastly, the interchange law holds as in (\ref{interchange1}) and
(\ref{interchange2}). For a ``one-sort formulation'' of pseudo
double category and mention of the subtleties in the following
remark, see \cite{grandisdouble1}.
\end{defn}

\begin{rmk}\label{weakeningremark}
The weakening of the horizontal 2-category to a horizontal
bicategory forces other parts of the notion of double category to
weaken in a pseudo double category $\mathbb{D}$. For example, the
horizontal composition of squares cannot be strictly associative if
the composition of horizontal morphisms is not strictly associative.
Similarly, the horizontal composition of squares cannot be strictly
unital if the composition of horizontal morphisms is not strictly
unital.

If the composition of horizontal morphisms is not strictly unital,
then $\mathbf{V}\mathbb{D}$ is neither a 2-category nor a
bicategory: the vertical composition of 2-cells in
$\mathbf{V}\mathbb{D}$ (which is the horizontal composition of
squares in $\mathbb{D}$) is not closed. However, if we redefine the
vertical composition of 2-cells $\beta$ and $\gamma$ in
$\mathbf{V}\mathbb{D}$ to be
$$\xymatrix@R=3pc@C=3pc{ \ar[rr]^{1^h} \ar[d] \ar@{}[rrd]|{\rho_{1^h}^{-1}} & & \ar[d] \\
\ar[r]^{1^h} \ar[d] \ar@{}[dr]|{\beta} &  \ar[r]^{1^h} \ar[d]
\ar@{}[dr]|{\gamma} & \ar[d] \\ \ar[r]_{1^h} \ar[d] \ar@{}[rrd]|{\rho_{1^h}} & \ar[r]_{1^h} & \ar[d] \\
\ar[rr]_{1^h} & & }$$ then we obtain a 2-category.

If $\lambda$ and $\rho$ are the vertical identity squares (\ie the
composition of horizontal morphisms is strictly unital), then
$i^h_j$ is a horizontal identity square by the naturality of
$\lambda$ and $\rho$, and hence $\mathbf{V}\mathbb{D}$ is a
2-category without any alterations. If additionally $\alpha$ is a
vertical identity, then we obtain the usual notion of double
category. Whenever $\lambda$ and $\rho$ are the vertical identity
squares, we say that $\mathbb{D}$ {\it has strict units}. Note that
for pseudo double categories we must require
$i^h_{1^v_A}=i^v_{1_A^h}$ even though this equality follows from the
other axioms in the case of strict double categories.
\end{rmk}

\begin{examp} \label{Rngs}
In the pseudo double category $\mathbb{R}\text{ng}$, objects are
commutative rings, horizontal morphisms from $A$ to $B$ are
$(B,A)$-bimodules, vertical morphisms are ring homomorphisms, while
squares $\beta$ with boundary as in (\ref{sourcetarget}) are group
homomorphisms $\xymatrix@1{\beta:f \ar[r] & g}$ such that
$\beta(bxa)=k(b)\beta(x)j(a)$ for all $b \in B, x \in f,$ and $a \in
A$. Composition of horizontal morphisms is tensor product of
bimodules, while composition of vertical morphisms is ordinary
function composition.
\end{examp}

\begin{defn}
A {\it pseudo double functor} $\xymatrix@1{F: \mathbb{D} \ar[r] &
\mathbb{E}}$ consists of maps
$$\xymatrix{Obj \hspace{1mm} \mathbb{D} \ar[r] & Obj \hspace{1mm}  \mathbb{E}}$$
$$\xymatrix{Hor \hspace{1mm} \mathbb{D} \ar[r] & Hor \hspace{1mm} \mathbb{E}}$$
$$\xymatrix{Ver \hspace{1mm} \mathbb{D} \ar[r] & Ver \hspace{1mm} \mathbb{E}}$$
$$\xymatrix{Squares \hspace{1mm} \mathbb{D} \ar[r] & Squares \hspace{1mm} \mathbb{E}}$$
which preserve all sources and targets and are compatible with
compositions and units in the following sense. The restriction
$\mathbf{H}F$ to the horizontal bicategory is a
homomorphism\footnote{A {\it homomorphism of bicategories} preserves
composition and units up to coherence isos. The composition
coherence iso satisfies a hexagon diagram with the composition
coherence isos of the bicategories, and the unit coherence iso
satisfies two square diagrams with the unit coherence isos of the
bicateogries. Some authors call this a {\it pseudo functor.}} of
bicategories whose coherence isos are natural with respect to {\it
all} squares, and further, the restriction of $F$ to the vertical
1-category is an ordinary functor. Equivalently, $F$ consists of
functors $\xymatrix@1{F_0:\mathbb{D}_0 \ar[r] & \mathbb{E}_0}$ and
$\xymatrix@1{F_1:\mathbb{D}_1 \ar[r] & \mathbb{E}_1}$ and natural
isomorphisms
$$\xymatrix@C=3pc@R=3pc{\mathbb{D}_1 \times_{\mathbb{D}_0}
\mathbb{D}_1 \ar[r] \ar[d]_{F_1
\times_{\mathbb{D}_0} F_1} & \mathbb{D}_1 \ar[d]^{F_1} \\
\mathbb{E}_1 \times_{\mathbb{E}_0} \mathbb{E}_1 \ar[r]
\ar@{=>}[ur]^{\gamma}
 & \mathbb{E}_1} \hspace{2cm} \xymatrix@C=3pc@R=3pc{\mathbb{D}_0 \ar[r]^{\eta} \ar[d]_{F_0}
 & \mathbb{D}_1 \ar[d]^{F_1} \\ \mathbb{E}_0 \ar[r]_{\eta} \ar@{=>}[ur]^{\delta} & \mathbb{E}_1}$$
whose components are squares with trivial vertical edges, and which
satisfy the usual three coherence diagrams for homomorphisms of
bicategories. The naturality of $\delta$ means
$$\begin{array}{c}
\xymatrix@C=3pc@R=3pc{FA \ar[r]^{1_{FA}^h} \ar@{=}[d]
\ar@{}[dr]|{\delta_A}
& FA \ar@{=}[d] \\ FA \ar[r]|{\lr{F1^h_A}} \ar[d]_{Fj} \ar@{}[dr]|{F(i^h_j)} & FA \ar[d]^{Fj} \\
FC \ar[r]_{F1_C^h} & FC}
\end{array}
=
\begin{array}{c}
\xymatrix@C=3pc@R=3pc{FA \ar[r]^{1_{FA}^h} \ar[d]_{Fj}
\ar@{}[dr]|{i^h_{Fj}} & FA \ar[d]^{Fj}
\\ FC \ar[r]|{\lr{1^h_{FC}}} \ar@{=}[d] \ar@{}[dr]|{\delta_C} & FC \ar@{=}[d]
\\ FC \ar[r]_{F1^h_C} & FC}
\end{array}
$$
for all vertical morphisms $j$. Thus if $\delta$ is trivial, then
$F$ preserves horizontal identity squares. If $\gamma$ is
additionally trivial and $\mathbb{D}$ and $\mathbb{E}$ are strict,
then $F$ is an internal functor in $Cat$, in other words a double
functor.
\end{defn}

\begin{defn}
Let $\mathbb{D}$ be a pseudo double category. A {\it pseudo
holonomy} is a homomorphism of bicategories
$\xymatrix@1{(\mathbf{V}\mathbb{D})_0 \ar[r] &
\mathbf{H}\mathbb{D}}$ which is the identity on objects.
\end{defn}

\begin{defn}
Let $\mathbb{D}$ be a pseudo double category. A {\it pseudo folding
on $\mathbb{D}$} consists of a pseudo holonomy $\xymatrix@1{j
\ar@{|->}[r] & \overline{j}}$ and bijections $\Lambda^{f,k}_{j,g}$
from squares in $\mathbb{D}$ with boundary as in (\ref{boundary1})
to squares in $\mathbb{D}$ with boundary as in (\ref{boundary2})
such that (i),(ii), (iii), and (iv) of Definition
\ref{foldingstructure} hold after composing with the coherence iso
2-cells of the horizontal bicategory and the pseudo holonomy. If the
pseudo holonomy of a pseudo folding is a strict 2-functor, then we
say simply {\it folding} instead of pseudo folding.
\end{defn}

\begin{rmk}
One would like to say that a pseudo folding on a pseudo double
category is a pseudo double functor $\xymatrix@1{\mathbb{D} \ar[r] &
\mathbb{Q}\mathbf{H}\mathbb{D}}$, but the quintets of a bicategory
unfortunately do not form a pseudo double category. Instead we write
out what the pseudo functoriality would mean directly: the holonomy
2-functor is replaced by a homomorphism of bicategories, and
composites of squares are preserved after composing with coherence
isos.
\end{rmk}

\begin{defn}
Let $\mathbb{D}$ and $\mathbb{E}$ be pseudo double categories with
pseudo folding. A {\it morphism of pseudo double categories with
pseudo folding} $\xymatrix@1{F:\mathbb{D} \ar[r] & \mathbb{E}}$ is a
morphism of pseudo double categories equipped with a coherence iso
2-cell
$$\xymatrix{F(\overline{j}) \ar[r] & \overline{F(j)}}$$
of the horizontal bicategory for each vertical morphism $j$ of
$\mathbb{D}$ such that:
\begin{enumerate}
\item
the coherence iso 2-cells are compatible with the coherence iso 2-cells of
the pseudo holonomies,
\item
after composing with the relevant coherence iso 2-cells, we have
$$F(\Lambda^{\mathbb{D}}(\alpha))=\Lambda^{\mathbb{E}}(F(\alpha))$$ for all
squares $\alpha$ of $\mathbb{D}$.
\end{enumerate}
\end{defn}

\begin{examp} \label{pseudofoldingexample}
Consider the pseudo double category $\mathbb{R}$ng of rings,
bimodules, ring homomorphisms, and squares in Example \ref{Rngs}.
From a ring homomorphism $\xymatrix@1{j:A \ar[r] & C}$ we get a
$(C,A)$-bimodule $\overline{j}$ by viewing $C$ as a left $C$-module
in the usual way and as a right $A$-module via $j$. We also denote
this bimodule by $C_j$ as in \cite{maysig}, where such base changes
are organized into a so-called closed symmetric bicategory.  The map
$\xymatrix@1{j \ar@{|->}[r] & \overline{j}}$ defines a pseudo
holonomy which strictly preserves units, but preserves compositions
only up to a coherence iso 2-cell. For a square
$$\xymatrix{A \ar[r]^M \ar[d]_j \ar@{}[dr]|\beta & B \ar[d]^k \\ C \ar[r]_N & D},$$
\ie for a group homomorphism $\xymatrix@1{\beta:M \ar[r] & N}$ such
that $\beta(bma)=k(b)\beta(m)j(a)$, define a 2-cell of
$\mathbf{H}\mathbb{R}$ng
$$\xymatrix{\Lambda(\beta):D_k \otimes_B M \ar@{=>}[r] & N \otimes_C C_j }$$
$$\xymatrix{ d \otimes m \ar@{|->}[r] & (d\cdot \beta(m)) \otimes 1_C.}$$
This is well defined because $\xymatrix@1{(d,m)\ar@{|->}[r] & d\cdot
\beta(m)}$ is middle $B$-linear.
$$\aligned
(d\cdot b)\cdot \beta(m) &= (dk(b))\cdot \beta(m) \\
&=d\cdot(k(b) \cdot \beta(m)) \\
&=d \cdot \beta(b\cdot m)
\endaligned$$
Then $\Lambda$ is bijective, a pre-image of a 2-cell
$\xymatrix@1{\beta':D_k \otimes_B M \ar@{=>}[r] & N \otimes_C C_j }$
is given by $\beta(m):=\beta'(1_D \otimes m)$ under the
identification $N\otimes_CC_j \cong N_j$. The coherence diagrams
associated with (i),(ii), and (iii) of Definition
\ref{foldingstructure} can be verified, and $\Lambda$ is a pseudo
folding on $\mathbb{R}$ng with pseudo holonomy $\xymatrix@1{j
\ar@{|->}[r] & \overline{j}}$. For more on this example in the
context of so-called anchored bicategories, trace maps, and
symmetric bicategories, see \cite{maysig}, \cite{pontotraces}, and
\cite{shulmananchor}.
\end{examp}

\begin{examp} \label{worldsheetexample}
The pseudo double category $\mathbb{W}$ of worldsheets is relevant
to conformal field theory, and admits a folding.  A {\it worldsheet}
$x$ is a real, compact, not necessarily connected, two dimensional,
smooth manifold with complex structure and real analytically
parametrized boundary components. A boundary component $k$ is called
{\it inbound} if the orientation of its parametrization
$\xymatrix@1{f_k:S^1 \ar[r] & k}$ with respect to the orientation on
$k$ is the same as the orientation of the identity parametrization
of the boundary of the unit disk. Otherwise $k$ is called {\it
outbound}. We say that the inbound components of $x$ are {\it
labelled by a finite set $A$} if $x$ is equipped with a bijection
between the set of inbound components and $A$.

The objects of $\mathbb{W}$ are finite sets, the horizontal
morphisms from $A$ to $B$ are worldsheets $x$ whose inbound
respectively outbound components are labelled by $A$ respectively
$B$, the vertical morphisms are bijections of sets. For two finite
sets $A$ and $B$ of the same cardinality, we also include in our
horizontal morphisms from $A$ to $B$ unions of unparametrized
circles $(a,S^1,b)$ where each $a\in A$ and each $b \in B$ appear
only once, so we can view such unions as bijections. The circle
$(a,S^1,b)$ is viewed as an infinitely thin annulus with inbound
component labelled by $a$ and outbound component labelled by $b$.
For $x$ and $y$ worldsheets, a square
$$\xymatrix{A \ar[r]^x \ar[d]_j \ar@{}[dr]|{\beta} & B \ar[d]^k \\ C \ar[r]_y & D }$$
consists of a holomorphic diffeomorphism $\xymatrix@1{\beta:x \ar[r]
& y}$ which:
\begin{enumerate}
\item
takes every inbound component of $x$ labelled by $a \in A$ to an
inbound component of $y$ labelled by $j(a)$,
\item
takes every outbound component of $x$ labelled by $b \in B$ to an
outbound component of $y$ labelled by $k(b)$, and
\item
preserves the boundary parametrizations.
\end{enumerate}
For $x$ and $y$ unions of circles, there is a unique square $\beta$
with boundary as above if and only if $x=y$, and $j$ and $k$ are
vertical identities. This unique square is necessarily the identity.
There are no other squares in this double category.

The composition of horizontal morphisms is given by gluing of
surfaces, and the horizontal composition of squares is defined
analogously. The composition of vertical morphisms is composition of
functions. The strict horizontal unit from $A$ to $A$ is the union
of circles $(a,S^1,a)$ for $a \in A$. A strict holonomy is defined
by mapping a function $j$ to the union of circles $(a,S^1,j(a))$.

Actually, this example is an illustration of a pseudo version of
Remark \ref{holonomy=mixedcomposition} and Remark \ref{extendingD}.
The worldsheets admit a mixed composition with the bijections via
relabelling. Including the disjoint union of circles also as
horizontal morphisms corresponds to extending $\mathbb{D}$ to
$\mathbb{D}'$, and this explains our choice of squares. The holonomy
is the inclusion.

A folding $\Lambda$ is given by simply relabelling the outbound
boundary components of $x$ via $k$ and relabelling the inbound
boundary components of $y$ via $j$: the holomorphic diffeomorphism
stays the same.

This example, along with \cite{kerlerlyubashenko2001} and
\cite{moskaliukvlassov1998}, suggests that double categories play a
role in the mathematics relating to field theories and high energy
physics.
\end{examp}

\section{Pseudo Algebras and Pseudo Double Categories with Folding}\label{pseudosection}

As we have seen in Theorems \ref{YZ}, \ref{XY}, and \ref{XZ}, the
2-categories of $I$-categories, double categories with folding, and
certain strict 2-functors are 2-equivalent if $I$ is a groupoid, and
the latter two remain 2-equivalent even if $I$ is merely a category.
Next we work towards pseudo versions of these theorems, as stated in
Theorems \ref{YZpseudo}, \ref{XYpseudo}, and \ref{XZpseudo}. We
prove the 2-equivalence of three 2-categories
$\mathcal{X},\mathcal{Y},$ and $\mathcal{Z}$ as
introduced\footnote{For the 2-equivalences with $\mathcal{X}$ we
assume $I$ is a groupoid.} in Notations \ref{Xpseudonotation},
\ref{Ypseudonotation}, and \ref{Zpseudonotation}. The 2-category
$\mathcal{X}$ is the 2-category of pseudo $I$-categories as defined
below, while $\mathcal{Y}$ and $\mathcal{Z}$ are the 2-categories of
certain pseudo double categories with pseudo folding and certain
2-functors respectively.

\begin{defn}
A {\it pseudo algebra over the 2-theory of categories with underlying category $I$}, also called
a {\it pseudo $I$-category} for short, is a category $I$ and a
strict 2-functor $\xymatrix@1{X:I^2 \ar[r] & Cat}$ with strictly
2-natural functors
$$\xymatrix@1{X_{B,C} \times X_{A,B} \ar[r]^-{\circ} & X_{A,C}}$$
$$\xymatrix@1{ \ast \ar[r]^{\eta_B} & X_{B,B} }$$
for all $A,B,C \in I$ and natural isomorphisms
$$\xymatrix@C=5pc{X_{C,D} \times (X_{B,C} \times X_{A,B}) \ar[r]^-{\circ \times 1_{X_{A,B}}}
\ar[dd]_{\cong} \ar@{}[ddr]|{\Downarrow \alpha_{A,B,C,D}} & X_{B,D}
\times
X_{A,B} \ar[dr]^{\circ} & \\
& &  X_{A,D} \\ (X_{C,D} \times X_{B,C}) \times X_{A,B}
\ar[r]_-{1_{X_{C,D}} \times \circ} & X_{C,D} \times X_{A,C}
\ar[ur]_{\circ} & }$$

\begingroup
\vspace{-2\abovedisplayskip} \Small
$$\begin{array}{cc}
\xymatrix@C=4pc@R=3pc{\ast \times X_{A,B} \ar[r]^-{\eta_B \times
1_{X_{A,B}}} \ar[d]_{pr_2} \ar@{}[dr]|{\Downarrow \lambda_{A,B}} &
X_{B,B} \times X_{A,B} \ar[d]^-{\circ} \\ X_{A,B}
\ar[r]_{1_{X_{A,B}}} & X_{A,B} } & \xymatrix@C=4pc@R=3pc{ X_{B,C}
\times \ast \ar@{}[dr]|{\Downarrow \rho_{B,C}} \ar[r]^-{1_{X_{B,C}}
\times \eta_B} \ar[d]_{pr_1} & X_{B,C} \times X_{B,B}
\ar[d]^-{\circ}
\\  X_{B,C} \ar[r]_{1_{X_{B,C}}} & X_{B,C} }
\end{array}
$$
\endgroup
\noindent which are the components of modifications and satisfy the
usual coherence diagrams for a bicategory as in the original
\cite{benabou1}, or in the review \cite{leinsterbicat}, or in the
Appendix to \cite{fiore3}. The requirement that $\alpha, \lambda,$
and $\rho$ be modifications means
$X_{j,m}(\alpha_{h,g,f})=\alpha_{X_{\ell,m}(h),X_{k,\ell}(g),X_{j,k}(f)}$,
$X_{j,k}(\lambda_f)=\lambda_{X_{j,k}(f)}$, and
$X_{k,\ell}(\rho_{g})=\rho_{X_{k,\ell}(g)}$ for all
$$\xymatrix{j:A \ar[r] & A'}$$
$$\xymatrix{k:B \ar[r] & B'}$$ $$\xymatrix{\ell:C \ar[r] & C'}$$ $$\xymatrix{m:D \ar[r] & D'.}$$
If $\lambda$ and $\rho$ are identities, we say $X$ {\it has strict
units}. We denote the value of $\eta_B$ on the unique object and
morphism of the terminal category by $1_B$ and $i_{1_B}$
respectively. We denote the identity morphism on an object $f$ in
the category $X_{A,B}$ by $i_f$.
\end{defn}

\begin{defn}
A {\it morphism  of pseudo $I$-categories} $\xymatrix@1{F: X \ar[r]
& Y}$ is a strict 2-natural transformation $\xymatrix@1{F:X
\ar@{=>}[r] & Y}$ with natural isomorphisms
$$\xymatrix@R=3pc@C=3pc{ X_{B,C} \times X_{A,B} \ar@{}[dr]|{\quad
\Uparrow \gamma_{A,B,C}} \ar[r]^-{\circ} \ar[d]_{F_{B,C} \, \times
\, F_{A,B}} & X_{A,C} \ar[d]^{F_{A,C}}
\\ Y_{B,C} \times Y_{A,B} \ar[r]_-{\circ} & Y_{A,C}}$$
$$\xymatrix@R=3pc@C=3pc{\ast \ar[r]^{\eta^X_A} \ar[d] \ar@{}[dr]|{\Uparrow \delta_A} & X_{A,A} \ar[d]^{F_{A,A}} \\
\ast \ar[r]_{\eta^Y_A} & Y_{A,A} }$$ which are the components of
modifications and satisfy the usual coherence diagrams for
homomorphisms of bicategories. The requirement that $\gamma$ and
$\delta$ be modifications is equivalent to
$Y_{j,\ell}(\gamma_{g,f})=\gamma_{Y_{k,\ell}(g),Y_{j,k}(f)}$ and
$Y_{j,j}(\delta_A)=\delta_{A'}.$
\end{defn}

\begin{defn}
A {\it 2-cell $\xymatrix@1{\sigma:F \ar@{=>}[r] & G}$ in the
2-category of pseudo $I$-categories} is a modification
$\xymatrix@1{\sigma:F \ar@{~>}[r] & G}$ compatible with composition
and identity. More specifically, a 2-cell $\sigma$ consists of
natural transformations $\xymatrix@1{\sigma_{A,B}:F_{A,B}
\ar@{=>}[r] & G_{A,B}}$ for all $A,B \in I$ such that
$$Y_{j,k}(\sigma_{A,B}^f)=\sigma_{C,D}^{X_{j,k}(f)}$$
$$\gamma^G_{g,f} \odot (\sigma^g_{B,C} \circ \sigma^f_{A,B})=\sigma_{A,C}^{g \circ f} \odot \gamma^F_{g,f}$$
$$\sigma_{A,A}^{1_A} \odot \delta^F_A=\delta^G_A$$ for all $\xymatrix@1{(j,k):(A,B) \ar[r] &
(C,D)}$ in $I^2, f \in X_{A,B}$, $g \in X_{B,C}$, and all objects
$A$ of $I$. Here $\odot$ denotes the composition in the categories
$X_{A,B}$ and $\circ$ denotes the composition functor of pseudo
$I$-categories.
\end{defn}

Following the convention introduced in Section
\ref{strictcategories}, we use the term {\it pseudo $I$-category} to
abbreviate {\it pseudo algebra over the 2-theory of categories with
underlying category $I$}. The morphisms and 2-cells above are the
morphisms and 2-cells in the 2-category of pseudo algebras over the
2-theory of categories with underlying category $I$ as in
\cite{fiore1}, \cite{hu}, and \cite{hu1}. In this section $I$ will
denote a fixed category. Whenever we require $I$ to be a groupoid,
we will explicitly say so.

\begin{rmk}
If $I$ is a groupoid, units are strict, and $F$ and $G$ are
morphisms of pseudo $I$-categories that strictly preserve the units,
then Remark \ref{holonomyconstruction} holds. In particular, the
requirement $Y_{j,k}(\sigma_{A,B}^f)=\sigma_{C,D}^{X_{j,k}(f)}$ on a
2-cell $\sigma$ can replaced by $\sigma_{A,C}^{P(j)}=i_{P'(j)}$.
\end{rmk}

\begin{notation} \label{Xpseudonotation}
Let $\mathcal{X}$ denote the 2-category of pseudo $I$-categories
with strict units. The morphisms of $\mathcal{X}$ are morphisms of
pseudo $I$-categories which preserve the units strictly.
\end{notation}

\begin{notation} \label{Ypseudonotation}
Let $\mathcal{Y}$ denote the 2-category of pseudo double categories
$\mathbb{D}$ with strict units equipped with a folding (\ie the
holonomy is strict) and such that $(\mathbf{V}\mathbb{D})_0=I$. We
further require that the associativity coherence iso
$\alpha_{\overline{k},f,\overline{\ell}}$ is the identity for
vertical morphisms $k$ and $\ell$ and horizontal morphisms $f$ such
that $\overline{k} \circ f \circ \overline{\ell}$ exists.

A morphism in $\mathcal{Y}$ is a morphism $F$ of pseudo double
categories with folding which preserves the holonomy and units
strictly and is the identity on $(\mathbf{V}\mathbb{D})_0$. We
further require that $\gamma^F_{\overline{k},f}$ and
$\gamma^F_{f,\overline{\ell}}$ are identities.

A 2-cell $\xymatrix@1{\sigma:F \ar@{=>}[r] & G}$ in $\mathcal{Y}$ is
a vertical natural transformation that is compatible with folding
and has identity components. Less succinctly, a 2-cell assigns to
each pair $(A,B) \in I^2$ a natural transformation
$\xymatrix@1{\sigma_{A,B}:\mathbf{H}F_{A,B} \ar@{=>}[r] &
\mathbf{H}G_{A,B}}$ such that
$$\sigma^{\overline{j}}_{A,C}=i_{\overline{j}}^v$$
$$\begin{bmatrix}
\ [\ \sigma^f_{A,B} \ \ \sigma^g_{B,C} \ ] \  \\ \gamma^G_{g,f}
\end{bmatrix}=\begin{bmatrix}\gamma^F_{g,f} \vspace{2mm}
\\  \ \sigma^{[f \ g]}_{A,C} \ \end{bmatrix}$$
$$\sigma_{A,A}^{1^h_A}=i_{1_A^h}$$
for all vertical morphisms $j \in I(A,C)$, composable horizontal
morphisms $f$ and $g$, and all objects $A$.
\end{notation}

\begin{notation} \label{Zpseudonotation}
Another 2-category of interest is the 2-category $\mathcal{Z}$. An
object of $\mathcal{Z}$ is a strict 2-functor $\xymatrix@1{P:I
\ar[r] & \mathcal{C}}$ into a bicategory $\mathcal{C}$ with strict
units which is the identity on objects. We further require that the
associativity coherence iso $\alpha_{P(k), f,P(\ell)}$ of the
bicategory $\mathcal{C}$ is the identity for morphisms $k$ and
$\ell$ of $I$ such that $P(k) \circ f \circ P(\ell)$ exists in
$\mathcal{C}$. The object set of $\mathcal{C}$ is $Obj \hspace{1mm}
I$.

A morphism from $\xymatrix@1{P:I \ar[r] & \mathcal{C}}$ to
$\xymatrix@1{P':I \ar[r] & \mathcal{C}'}$ in $\mathcal{Z}$ is a
homomorphism of bicategories $\xymatrix@1{F:\mathcal{C} \ar[r] &
\mathcal{C}'}$ which strictly preserves units and such that
\begin{equation} \label{Zmorphism}
\xymatrix@C=4pc@R=3pc{I \ar[r]^P \ar[dr]_{P'} & \mathcal{C} \ar[d]^F
\\ & \mathcal{C'}}
\end{equation}
strictly commutes. We further require that $\gamma^F_{P(k),f}$ and
$\gamma^F_{f,P(\ell)}$ are identities.

A 2-cell $\xymatrix@1{\sigma:F \ar@{=>}[r] & G}$ consists of natural
transformations $\sigma_{A,B}$ for all $A,B \in Obj \hspace{1mm}
\mathcal{C}$ such that
$$\sigma_{A,C}^{P(j)}=i_{P'(j)}$$
$$\gamma^G_{g,f}\odot(\sigma^g_{B,C} \circ \sigma^f_{A,B})=\sigma^{g \circ f}_{A,C}\odot \gamma^F_{g,f}$$
$$\sigma_{A,A}^{1^h_A}=i_{1_A^h}$$
for all $j \in I(A,C), f \in \mathcal{C}(A,B), g \in
\mathcal{C}(B,C),$ and all objects $A$ of $I$. Here $\odot$ denotes
the vertical composition of 2-cells in a bicategory, while $\circ$
denotes the horizontal composition of 2-cells.
\end{notation}

\begin{rmk}
The requirement that units be strict in a pseudo double category is
not as rigid as it first seems, since this can be arranged in most
examples. The authors of \cite{grandisdouble1} and
\cite{grandisdouble2} also assume that units are strict, and arrange
it in most of their examples.
\end{rmk}

\begin{thm} \label{YZpseudo}
Let $I$ be a category. The 2-category $\mathcal{Y}$ of pseudo double
categories $\mathbb{D}$ with strict units equipped with a folding
such that $(\mathbf{V}\mathbb{D})_0=I$ is 2-equivalent to the
2-category $\mathcal{Z}$ of strict 2-functors $\xymatrix@1{I \ar[r]
& \mathcal{C}}$ into bicategories $\mathcal{C}$ with strict units
which are the identity on objects as in Notation
\ref{Ypseudonotation} and Notation \ref{Zpseudonotation}.
\end{thm}

\begin{thm} \label{XYpseudo}
Let $I$ be a groupoid. The 2-category $\mathcal{X}$ of pseudo
$I$-categories with strict units (pseudo algebras over the 2-theory
of categories with underlying groupoid $I$ and strict units) is
2-equivalent to the 2-category $\mathcal{Y}$ of pseudo double
categories $\mathbb{D}$ with strict units equipped with a folding
such that $(\mathbf{V}\mathbb{D})_0=I$ as defined in Notation
\ref{Xpseudonotation} and Notation \ref{Ypseudonotation}.
\end{thm}


\begin{thm} \label{XZpseudo}
Let $I$ be a groupoid. The 2-category $\mathcal{X}$ of pseudo
$I$-categories with strict units (pseudo algebras over the 2-theory
of categories with underlying groupoid $I$ and strict units) is
2-equivalent to the 2-category $\mathcal{Z}$ of strict 2-functors
$\xymatrix@1{I \ar[r] & \mathcal{C}}$ into bicategories
$\mathcal{C}$ with strict units which are the identity on objects as
in Notation \ref{Xpseudonotation} and Notation
\ref{Zpseudonotation}.
\end{thm}

We omit the proofs of Theorems \ref{YZpseudo}, \ref{XYpseudo}, and
\ref{XZpseudo} since they are straightforward but tedious
elaborations of the strict Theorems \ref{YZ}, \ref{XY}, and
\ref{XZ}. The strictness of units for $X$ in $\mathcal{X}$
corresponds to the strictness of the holonomy in $\mathcal{Y}$ and
the strictness of $\xymatrix@1{P:I \ar[r] & \mathcal{C}}$ in
$\mathcal{Z}$. The fact that morphisms of $\mathcal{X}$ strictly
preserve units corresponds to strict preservation of holonomy by
morphisms in $\mathcal{Y}$, as well as the strict preservation of
units by morphisms in $\mathcal{Z}$ and the strict commutativity of
Diagram (\ref{Zmorphism}).

\begin{examp} \label{ringsinY}
The pseudo double category $\mathbb{R}$ng in Example
\ref{pseudofoldingexample} can be slightly modified to make it into
an object of $\mathcal{Y}$ in Theorem \ref{YZpseudo} and Theorem
\ref{XYpseudo}. First we require the vertical morphisms to be
isomorphisms of rings, then note that bimodules admit a mixed
composition with isomorphisms of rings, and apply Remark
\ref{holonomy=mixedcomposition} and Remark \ref{extendingD}. Thus
the horizontal morphisms of $\mathbb{R}\text{ng}'_{\text{iso}}$ are
bimodules as well as isomorphisms of rings. A $(B,A)$-bimodule $M$
is composed with a ring isomorphism $\xymatrix@1{k:B \ar[r] & D}$ to
give a $(D,A)$-bimodule $k \circ M$ with underlying abelian group
$M$ by defining $d \cdot m:=k^{-1}(d)\cdot m$. The composition $N
\circ j$ is defined similarly. The squares of
$\mathbb{R}\text{ng}'_{\text{iso}}$ are the squares of
$\mathbb{R}$ng with invertible vertical sides, along with vertical
identities of the isomorphisms of rings. The holonomy is then an
inclusion and the horizontal bicategory is strictly unital.
\end{examp}

\begin{examp}
The pseudo double category $\mathbb{W}$ of worldsheets in Example
\ref{worldsheetexample} is an object of $\mathcal{Y}$ in Theorem
\ref{YZpseudo} and Theorem $\ref{XYpseudo}$ with $I$ the category of
finite sets and bijections. The horizontal morphisms are worldsheets
as well as bijections.
\end{examp}

\begin{thm}
Analogues of Theorems \ref{YZpseudo}, \ref{XYpseudo}, and
\ref{XZpseudo} hold for weak units and pseudo foldings, though
``2-equivalence'' must be replaced by ``biequivalence.'' Pseudo
$I$-categories with weak units correspond to pseudo double
categories with weak units and pseudo foldings, which in turn
correspond to homomorphisms of bicategories $P$ from the groupoid
$I$ to a bicategory $\mathcal{C}$ with weak units. Morphisms of
pseudo $I$-categories then correspond to morphisms of pseudo double
categories that preserve the pseudo holonomy up to coherence iso,
which in turn correspond to homomorphisms $F$ of bicategories such
that (\ref{Zmorphism}) commutes on objects strictly, but has a
coherence iso 2-cell $FP(j) \cong P'(j)$ for each morphism $j$ of
$I$.
\end{thm}

\begin{pf}
Omitted. The proof relies on a construction like $L(P)$ in the proof
of Theorem 6.5 in \cite{fiore3} to remedy
$$[ [ \overline{\ell} \ f  ] \ \overline{k}] \neq [ \overline{\ell}\ [ f
\ \overline{k} ] ]$$
$$F([ \overline{\ell} \ f \ \overline{k} ]) \neq [F(\overline{\ell}) \ F(f) \
F(\overline{k})]$$
$$(P(k) \circ f) \circ P(\ell) \neq P(k) \circ (f \circ P(\ell))$$
$$F(P(k) \circ f \circ P(\ell)) \neq F(P(k)) \circ F(f) \circ F(P(\ell)).$$
\end{pf}

This completes our comparison of strict 2-algebras and pseudo
algebras over the 2-theory of categories with variants of double
categories and 2-functors $\xymatrix@1{I \ar[r] & \mathcal{C}}$.

\def\cprime{$'$}


\begin{thebibliography}{10}

\bibitem{adamekrosicky1994}
Ji{\v{r}}{\'{\i}} Ad{\'a}mek and Ji{\v{r}}{\'\i} Rosick{\'y}.
\newblock {\em Locally presentable and accessible categories}, volume 189 of
  {\em London Mathematical Society Lecture Note Series}.
\newblock Cambridge University Press, Cambridge, 1994.

\bibitem{alaglthesis}
Fahd~Ali Al-Agl.
\newblock {\em Aspects of Multiple Categories}.
\newblock PhD thesis, University of Wales, 1989.

\bibitem{alglbrownsteiner2002}
Fahd~Ali Al-Agl, Ronald Brown, and Richard Steiner.
\newblock Multiple categories: the equivalence of a globular and a cubical
  approach.
\newblock {\em Adv. Math.}, 170(1):71--118, 2002.

\bibitem{andruskiewitschnatalequantum}
Nicol{\'a}s Andruskiewitsch and Sonia Natale.
\newblock Double categories and quantum groupoids.
\newblock {\em Publ. Mat. Urug.}, 10:11--51, 2005.

\bibitem{andruskiewitschnataletensor}
Nicol{\'a}s Andruskiewitsch and Sonia Natale.
\newblock Tensor categories attached to double groupoids.
\newblock {\em Adv. Math.}, 200(2):539--583, 2006.

\bibitem{baez}
John~C. Baez and James Dolan.
\newblock Categorification.
\newblock In {\em Higher category theory (Evanston, IL, 1997)}, volume 230 of
  {\em Contemp. Math.}, pages 1--36. Amer. Math. Soc., Providence, RI, 1998.

\bibitem{baezlauda}
John~C. Baez and Aaron~D. Lauda.
\newblock Higher-dimensional algebra. {V}. 2-groups.
\newblock {\em Theory Appl. Categ.}, 12:423--491 (electronic), 2004.

\bibitem{ehresmannone}
Andr{\'e}e Bastiani and Charles Ehresmann.
\newblock Multiple functors. {I}. {L}imits relative to double categories.
\newblock {\em Cahiers Topologie G\'eom. Diff\'erentielle}, 15(3):215--292,
  1974.

\bibitem{benabou1}
Jean B{\'e}nabou.
\newblock Introduction to bicategories.
\newblock In {\em Reports of the Midwest Category Seminar}, pages 1--77.
  Springer, Berlin, 1967.

\bibitem{benabou}
Jean B{\'e}nabou.
\newblock Structures alg\'ebriques dans les cat\'egories.
\newblock {\em Cahiers Topologie G\'eom. Diff\'erentielle}, 10:1--126, 1968.

\bibitem{bergnermultisorted}
Julia~E. Bergner.
\newblock Rigidification of algebras over multi-sorted theories.
\newblock {\em Preprint}, 2005.

\bibitem{boardmanvogt1973}
J.M. Boardman and R.M. Vogt.
\newblock {\em Homotopy invariant algebraic structures on topological spaces}.
\newblock Springer-Verlag, Berlin, 1973.
\newblock Lecture Notes in Mathematics, Vol. 347.

\bibitem{brown99}
Ronald Brown.
\newblock Groupoids and crossed objects in algebraic topology.
\newblock {\em Homology Homotopy Appl.}, 1:1--78 (electronic), 1999.

\bibitem{browngilbert1989}
Ronald Brown and N.~D. Gilbert.
\newblock Algebraic models of {$3$}-types and automorphism structures for
  crossed modules.
\newblock {\em Proc. London Math. Soc. (3)}, 59(1):51--73, 1989.

\bibitem{brownhigginscubes}
Ronald Brown and Philip~J. Higgins.
\newblock On the algebra of cubes.
\newblock {\em J. Pure Appl. Algebra}, 21(3):233--260, 1981.

\bibitem{brownhigginstensor}
Ronald Brown and Philip~J. Higgins.
\newblock Tensor products and homotopies for {$\omega$}-groupoids and crossed
  complexes.
\newblock {\em J. Pure Appl. Algebra}, 47(1):1--33, 1987.

\bibitem{brownicenhomotopies}
Ronald Brown and {\.I}lhan {\.I}{\c{c}}en.
\newblock Homotopies and automorphisms of crossed modules of groupoids.
\newblock {\em Appl. Categ. Structures}, 11(2):185--206, 2003.

\bibitem{brownicen2003}
Ronald Brown and {\.I}lhan {\.I}{\c{c}}en.
\newblock Towards a 2-dimensional notion of holonomy.
\newblock {\em Adv. Math.}, 178(1):141--175, 2003.

\bibitem{brownmackenzie}
Ronald Brown and Kirill~C.H. Mackenzie.
\newblock Determination of a double {L}ie groupoid by its core diagram.
\newblock {\em J. Pure Appl. Algebra}, 80(3):237--272, 1992.

\bibitem{brownmosa99}
Ronald Brown and Ghafar~H. Mosa.
\newblock Double categories, {$2$}-categories, thin structures and connections.
\newblock {\em Theory Appl. Categ.}, 5:No. 7, 163--175 (electronic), 1999.

\bibitem{brownspencer76}
Ronald Brown and Christopher~B. Spencer.
\newblock Double groupoids and crossed modules.
\newblock {\em Cahiers Topologie G\'eom. Diff\'erentielle}, 17(4):343--362,
  1976.

\bibitem{brownspencer74}
Ronald Brown and Christopher~B. Spencer.
\newblock {$G$}-groupoids, crossed modules and the fundamental groupoid of a
  topological group.
\newblock {\em Nederl. Akad. Wetensch. Proc. Ser. A (Indag. Math.)},
  38(4):296--302, 1976.

\bibitem{cheng1}
Eugenia Cheng and Aaron Lauda.
\newblock {\em Higher-Dimensional Categories: an illustrated guide book}.
\newblock Preprint.

\bibitem{datuashvili}
Tamar Datuashvili.
\newblock Whitehead homotopy equivalence and internal category equivalence of
  crossed modules in categories of groups with operations.
\newblock {\em Proc. A. Razmadze Math. Inst.}, 113:3--30, 1995.

\bibitem{dawsonparepronk2003adjoining}
R.J.~MacG. Dawson, R.~Par{\'e}, and D.A. Pronk.
\newblock Adjoining adjoints.
\newblock {\em Adv. Math.}, 178(1):99--140, 2003.

\bibitem{dawsonparepronkundecidable2003}
R.J.~MacG. Dawson, R.~Par{\'e}, and D.A. Pronk.
\newblock Undecidability of the free adjoint construction.
\newblock {\em Appl. Categ. Structures}, 11(5):403--419, 2003.

\bibitem{dawsonparepronkfree2004}
R.J.~MacG. Dawson, R.~Par{\'e}, and D.A. Pronk.
\newblock Free extensions of double categories.
\newblock {\em Cah. Topol. G\'eom. Diff\'er. Cat\'eg.}, 45(1):35--80, 2004.

\bibitem{dawsonparepronkspan2006}
R.J.~MacG. Dawson, R.~Par{\'e}, and D.A. Pronk.
\newblock More general spans.
\newblock {\em Preprint}, 2006.

\bibitem{dawsonparepronkpathology2006}
R.J.~MacG. Dawson, R.~Par{\'e}, and D.A. Pronk.
\newblock The pathology of double categories.
\newblock {\em Preprint}, 2006.

\bibitem{dawsonpare1993}
Robert Dawson and Robert Par{\'e}.
\newblock General associativity and general composition for double categories.
\newblock {\em Cahiers Topologie G\'eom. Diff\'erentielle Cat\'eg.},
  34(1):57--79, 1993.

\bibitem{dawsonpare2002free}
Robert Dawson and Robert Par{\'e}.
\newblock What is a free double category like?
\newblock {\em J. Pure Appl. Algebra}, 168(1):19--34, 2002.

\bibitem{ehresmanntwo}
Andr{\'e}e Ehresmann and Charles Ehresmann.
\newblock Multiple functors. {II}. {T}he monoidal closed category of multiple
  categories.
\newblock {\em Cahiers Topologie G\'eom. Diff\'erentielle}, 19(3):295--333,
  1978.

\bibitem{ehresmannthree}
Andr{\'e}e Ehresmann and Charles Ehresmann.
\newblock Multiple functors. {III}. {T}he {C}artesian closed category {${\rm
  Cat}\sb{n}$}.
\newblock {\em Cahiers Topologie G\'eom. Diff\'erentielle}, 19(4):387--443,
  1978.

\bibitem{ehresmannfour}
Andr{\'e}e Ehresmann and Charles Ehresmann.
\newblock Multiple functors. {IV}. {M}onoidal closed structures on {${\rm
  Cat}\sb{n}$}.
\newblock {\em Cahiers Topologie G\'eom. Diff\'erentielle}, 20(1):59--104,
  1979.

\bibitem{ehresmann}
Charles Ehresmann.
\newblock Cat\'egories structur\'ees.
\newblock {\em Ann. Sci. \'Ecole Norm. Sup. (3)}, 80:349--426, 1963.

\bibitem{ehresmannquintett63}
Charles Ehresmann.
\newblock Cat\'egories structur\'ees. {III}. {Q}uintettes et applications
  covariantes.
\newblock In {\em Topol. et G\'eom. Diff. (S\'em. C. Ehresmann), Vol. 5},
  page~21. Institut H. Poincar\'e, Paris, 1963.

\bibitem{ehresmann2}
Charles Ehresmann.
\newblock {\em Cat\'egories et structures}.
\newblock Dunod, Paris, 1965.

\bibitem{fiore1}
Thomas~M. Fiore.
\newblock Pseudo limits, biadjoints, and pseudo algebras: categorical
  foundations of conformal field theory.
\newblock {\em Mem. Amer. Math. Soc.}, 182(860), 2006,
  http://arxiv.org/abs/math.CT/0408298.

\bibitem{fiore3}
Thomas~M. Fiore.
\newblock On the cobordism and commutative monoid with cancellation approaches
  to conformal field theory.
\newblock {\em J. Pure Appl. Algebra}, 209(3):583--620, 2007.

\bibitem{fiorekrizhu}
Thomas~M. Fiore, Igor Kriz, and Po~Hu.
\newblock Pseudo algebras with {L}aplaza sets.
\newblock {\em Preprint}.

\bibitem{forrester}
Magnus Forrester-Barker.
\newblock Group objects and internal categories.
\newblock 2002, http://arxiv.org/abs/math.CT/0212065.

\bibitem{garnerthesis}
Richard Garner.
\newblock {\em Polycategories}.
\newblock PhD thesis, University of Cambridge.

\bibitem{garnerpseudodistributive2005}
Richard Garner.
\newblock Polycategories via pseudo-distributive laws.
\newblock http://www.dpmms.cam.ac.uk/{$\sim$}rhgg2/.

\bibitem{garner2005}
Richard Garner.
\newblock Double clubs.
\newblock {\em Cahiers Topologie G\'eom. Diff\'erentielle}, 47(4):261--317,
  2006.

\bibitem{grandisdouble1}
Marco Grandis and Robert Par\'e.
\newblock Limits in double categories.
\newblock {\em Cahiers Topologie G\'eom. Diff\'erentielle Cat\'eg.},
  40(3):162--220, 1999.

\bibitem{grandisdouble2}
Marco Grandis and Robert Par\'e.
\newblock Adjoints for double categories. {A}ddenda to: ``{L}imits in double
  categories''.
\newblock {\em Cah. Topol. G\'eom. Diff\'er. Cat\'eg.}, 45(3):193--240, 2004.

\bibitem{higgins1963}
Philip~J. Higgins.
\newblock Algebras with a scheme of operators.
\newblock {\em Math. Nachr.}, 27:115--132, 1963.

\bibitem{higgins2005}
Philip~J. Higgins.
\newblock Thin elements and commutative shells in cubical
  {$\omega$}-categories.
\newblock {\em Theory Appl. Categ.}, 14:No. 4, 60--74 (electronic), 2005.

\bibitem{hu}
Po~Hu and Igor Kriz.
\newblock Conformal field theory and elliptic cohomology.
\newblock {\em Adv. Math.}, 189(2):325--412, 2004,
  http://www.math.lsa.umich.edu/{$\sim$}ikriz/.

\bibitem{hu1}
Po~Hu and Igor Kriz.
\newblock Closed and open conformal field theories and their anomalies.
\newblock {\em Comm. Math. Phys.}, 254(1):221--253, 2005,
  http://www.math.lsa.umich.edu/{$\sim$}ikriz/.

\bibitem{kelly}
G.M. Kelly and Ross Street.
\newblock Review of the elements of {$2$}-categories.
\newblock In {\em Category Seminar (Proc. Sem., Sydney, 1972/1973)}, pages
  75--103. Lecture Notes in Math., Vol. 420. Springer, Berlin, 1974.

\bibitem{kerlerlyubashenko2001}
Thomas Kerler and Volodymyr~V. Lyubashenko.
\newblock {\em Non-semisimple topological quantum field theories for
  3-manifolds with corners}, volume 1765 of {\em Lecture Notes in Mathematics}.
\newblock Springer-Verlag, Berlin, 2001.

\bibitem{lawvere}
F.~William Lawvere.
\newblock Functorial semantics of algebraic theories.
\newblock {\em Proc. Nat. Acad. Sci. U.S.A.}, 50:869--872, 1963.

\bibitem{leinsterbicat}
Tom Leinster.
\newblock Basic bicategories.
\newblock 1998, http://arxiv.org/abs/math.CT/9810017.

\bibitem{leinster1}
Tom Leinster.
\newblock A survey of definitions of {$n$}-category.
\newblock {\em Theory Appl. Categ.}, 10:1--70 (electronic), 2002.

\bibitem{leinsteroperads2004}
Tom Leinster.
\newblock {\em Higher operads, higher categories}, volume 298 of {\em London
  Mathematical Society Lecture Note Series}.
\newblock Cambridge University Press, Cambridge, 2004,
  http://arxiv.org/abs/math.CT/0305049.

\bibitem{lodayfinitelymany}
Jean-Louis Loday.
\newblock Spaces with finitely many nontrivial homotopy groups.
\newblock {\em J. Pure Appl. Algebra}, 24(2):179--202, 1982.

\bibitem{maclane1}
Saunders Mac~Lane.
\newblock Natural associativity and commutativity.
\newblock {\em Rice Univ. Studies}, 49(4):28--46, 1963.

\bibitem{maclane1991}
Saunders Mac~Lane.
\newblock Coherence theorems and conformal field theory.
\newblock In {\em Category theory 1991 (Montreal, PQ, 1991)}, volume~13 of {\em
  CMS Conf. Proc.}, pages 321--328. Amer. Math. Soc., Providence, RI, 1992.

\bibitem{maclane3}
Saunders Mac~Lane.
\newblock {\em Categories for the working mathematician}, volume~5 of {\em
  Graduate Texts in Mathematics}.
\newblock Springer-Verlag, New York, second edition, 1998.

\bibitem{maclanewhitehead50}
Saunders Mac~Lane and J.H.C. Whitehead.
\newblock On the {$3$}-type of a complex.
\newblock {\em Proc. Nat. Acad. Sci. U. S. A.}, 36:41--48, 1950.

\bibitem{martins2004}
N.~Martins-Ferreira.
\newblock Weak categories in additive 2-categories with kernels.
\newblock In {\em Galois theory, Hopf algebras, and semiabelian categories},
  volume~43 of {\em Fields Inst. Commun.}, pages 387--410. Amer. Math. Soc.,
  Providence, RI, 2004.

\bibitem{martinspseudocategories}
N.~Martins-Ferreira.
\newblock Pseudo-categories.
\newblock {\em Journal of Homotopy and Related Structures}, 1(1):47--78, 2006,
  http://jhrs.rmi.acnet.ge/.

\bibitem{maysig}
J.P. May and J.\ Sigurdsson.
\newblock {\em Parametrized Homotopy Theory}, volume 132 of {\em Mathematical
  Surveys and Monographs}.
\newblock American Mathematical Society, Providence, RI, first edition, 2006.

\bibitem{mortondouble}
Jeffrey Morton.
\newblock A double bicategory of cobordisms with corners.
\newblock http://arxiv.org/abs/math.CT/0611930.

\bibitem{moskaliukvlassov1998}
S.~S. Moskaliuk and A.~T. Vlassov.
\newblock Double categories in mathematical physics.
\newblock {\em Ukra\"\i n. F\=\i z. Zh.}, 43(6-7):836--841, 1998.
\newblock International Symposium on Mathematical and Theoretical Physics
  (Kyiv, 1997).

\bibitem{palmquistthesis}
Paul~H. Palmquist.
\newblock {\em The Double Category of Adjoint Squares}.
\newblock PhD thesis, University of Chicago, 1969.

\bibitem{palmquist1971}
Paul~H. Palmquist.
\newblock The double category of adjoint squares.
\newblock In {\em Reports of the Midwest Category Seminar, V (Z\"urich, 1970)},
  Lecture Notes in Mathematics, Vol. 195, pages 123--153. Springer, Berlin,
  1971.

\bibitem{paolisurvey}
Simona Paoli.
\newblock Internal categorical structures in homotopical algebra.
\newblock In {\em Proceedings of the IMA 2004 Summer Program ``$n$-Categories:
  Foundations and Applications'' at the University of Minnesota}, To Appear.

\bibitem{pontotraces}
Kathleen Ponto.
\newblock {\em Fixed Point Theory and Trace for Bicategories}.
\newblock PhD thesis, University of Chicago, 2007.

\bibitem{segal1}
Graeme Segal.
\newblock The definition of conformal field theory.
\newblock In {\em Topology, geometry and quantum field theory}, volume 308 of
  {\em London Math. Soc. Lecture Note Ser.}, pages 421--577. Cambridge Univ.
  Press, Cambridge, 2004.

\bibitem{shulmananchor}
Michael Shulman.
\newblock {\em Framed and anchored bicategories}.
\newblock Preprint. 2007.

\bibitem{spencer77}
Christopher~B. Spencer.
\newblock An abstract setting for homotopy pushouts and pullbacks.
\newblock {\em Cahiers Topologie G\'eom. Diff\'erentielle}, 18(4):409--429,
  1977.

\bibitem{streetstructures}
Ross Street.
\newblock Categorical structures.
\newblock In {\em Handbook of algebra, Vol.\ 1}, pages 529--577. North-Holland,
  Amsterdam, 1996.

\bibitem{veritythesis}
Dominic Verity.
\newblock {\em Enriched Categories, Internal Categories and Change of Base}.
\newblock PhD thesis, University of Cambridge, 1992.

\bibitem{whitehead49}
J.H.C.\ Whitehead.
\newblock Combinatorial homotopy. {II}.
\newblock {\em Bull. Amer. Math. Soc.}, 55:453--496, 1949.

\end{thebibliography}
\end{document}